\ifx\LocalElevenPtArticleFormatLoaded\undefined
  \documentclass[a4paper,11pt]{article}     
  \usepackage{array}                        
  \usepackage{theorem}                      
  \usepackage{amsmath,amscd,amssymb}        
  \usepackage{latexsym}                     
  \usepackage[german,english]{babel}        
  \usepackage{graphicx}                     
  \usepackage{epsfig}
  \usepackage{xypic}
\fi

\usepackage{theorem}
\newtheorem{theorem}{Theorem}[section]
\newtheorem{lemma}[theorem]{Lemma}
\newtheorem{proposition}[theorem]{Proposition}
\newtheorem{cor}[theorem]{Corollary}
\theorembodyfont{\rm}

\renewcommand{\thetheorem}{\Roman{section}.\arabic{subsection}.\arabic{theorem}}

\newtheorem{definition}{Definition}

\newtheorem{problem}{Problem}        
\newtheorem{conjecture}{Conjecture}

\newcommand{\CC}{{\mathbb{C}}}

\newcommand{\HH}{{\mathbb{H}}}
\newcommand{\FF}{{\mathbb{F}}}
\newcommand{\NN}{{\mathbb{N}}}
\newcommand{\PP}{{\mathbb{P}}}
\newcommand{\QQ}{{\mathbb{Q}}}
\newcommand{\RR}{{\mathbb{R}}}
\newcommand{\ZZ}{{\mathbb{Z}}}

\newcommand{\calA}{{\cal A}}

\newcommand{\calL}{{\cal L}}

\newcommand{\cO}{{\cal O}}
\newcommand{\cA}{{\cal A}}
\newcommand{\cN}{{\cal N}}
\newcommand{\cX}{{\cal X}}
\newcommand{\on}[1]{\operatorname{#1}}
\newcommand{\Proj}{\on{Proj}}
\newcommand{\Spec}{\on{Spec}}
\newcommand{\PSp}{\on{PSp}}
\newcommand{\Sym}{\on{Sym}}
\newcommand{\BKm}{{\widetilde{\mathrm{Km}}\,}}
\newcommand{\cB}{{\cal B}}

\renewcommand{\Sp}{\on{Sp}}

\newcommand{\Ker}{\on{ker}}
\newcommand{\cI}{{\cal I}}
\newcommand{\To}{{\longrightarrow}}
\newcommand{\cM}{{\cal M}}
\newcommand{\cP}{{\cal P}}
\newcommand{\cS}{{\cal S}}

\newcommand\ratmap{\mathrel{{\hbox{\kern2pt\vrule height2.45pt depth-2.15pt
 width2pt}\kern1pt {\vrule height2.45pt depth-2.15pt width2pt}
  \kern1pt{\vrule height2.45pt depth-2.15pt width1.7pt\kern-1.7pt}
   {\raise1.4pt\hbox{$\scriptscriptstyle\succ$}}\kern1pt}}}
\newcommand{\SL}{\on{SL}}
\newcommand{\GL}{\on{GL}}
\newcommand{\PSL}{\on{PSL}}
\newcommand{\Gr}{\on{Gr}}
\renewcommand{\Im}{\on{Im}}
\newcommand{\eps}{\varepsilon}
\newcommand{\lev}{{\mbox{\tiny{\rm lev}}}}
\newcommand{\opn}{{\vphantom{\cA}^{\circ}\!\!}}
\font\gothic=eufm10
\newcommand \gothm {{\mbox{\gothic m}}}
\newenvironment{Proof}{\begin{ProofwCaption}{Proof}}{\end{ProofwCaption}}
\newenvironment{Proof*}[1]{\begin{ProofwCaption}{{#1}}}{\end{ProofwCaption}}
\newenvironment{ProofwCaption}[1]%
  {\addvspace\theorempreskipamount \noindent{\it #1.}\rm}%
  {\qed \par \addvspace\theorempostskipamount}
\newcommand{\qedsymbol}{\mbox{$\Box$}}
\newcommand{\qed}{\hfill\quad\qedsymbol}
\begin{document}
\title{The Geometry of Siegel Modular Varieties}
\date{}
\author{K.~Hulek and G.K.~Sankaran}
\maketitle
\begin{enumerate}
\item[] Introduction\ \dotfill\ \pageref{intro}
\item[I] Siegel modular varieties\ \dotfill\ \pageref{I}
 \begin{enumerate}
  \item[I.1] Arithmetic quotients of the Siegel upper half
plane\ \dotfill\ \pageref{I.1}
  \item[I.2] Compactifications of Siegel modular
varieties\ \dotfill\ \pageref{I.2}
 \end{enumerate}
\item[II] Classification theory\ \dotfill\ \pageref{II}
 \begin{enumerate}
   \item[II.1] The canonical divisor\ \dotfill\ \pageref{II.1}
   \item[II.2] The Kodaira dimension of
${\cal A}_g(n)$\ \dotfill\ \pageref{II.2}
   \item[II.3] Fundamental groups\ \dotfill\ \pageref{II.3}
 \end{enumerate}
\item[III] Abelian surfaces\ \dotfill\ \pageref{III}
 \begin{enumerate}
   \item[III.1] The lifting method\ \dotfill\ \pageref{III.1}
   \item[III.2] General type results for moduli spaces of abelian
surfaces\ \dotfill\ \pageref{III.2}
   \item[III.3] Left and right neighbours\ \dotfill\ \pageref{III.3}
 \end{enumerate}
\item[IV] Projective models\ \dotfill\ \pageref{IV}
  \begin{enumerate}
   \item[IV.1] The Segre cubic\ \dotfill\ \pageref{IV.1}
   \item[IV.2] The Burkhardt quartic\ \dotfill\ \pageref{IV.2}
   \item[IV.3] The Nieto quintic\ \dotfill\ \pageref{IV.3}
 \end{enumerate}
\item[V] Non-principal polarizations\ \dotfill\ \pageref{V}
 \begin{enumerate}
   \item[V.1] Type $(1,5)$ and the Horrocks-Mumford
bundle\ \dotfill\ \pageref{V.1}
   \item[V.2] Type $(1,7)$\ \dotfill\ \pageref{V.2}
   \item[V.3] Type $(1,11)$\ \dotfill\ \pageref{V.3}
   \item[V.4] Other type $(1,t)$ cases\ \dotfill\ \pageref{V.4}
 \end{enumerate}
\item[VI] Degenerations\ \dotfill\ \pageref{VI}
 \begin{enumerate}
   \item[VI.1] Local degenerations\ \dotfill\ \pageref{VI.1}
   \item[VI.2] Global degenerations and
compactification\ \dotfill\ \pageref{VI.2}
 \end{enumerate}
\item[] References\ \dotfill\ \pageref{biblio}
\end{enumerate}

\section{Introduction}\label{intro}

Siegel modular varieties are interesting because they arise as moduli
spaces for abelian varieties with a polarization and a level structure, and
also because of their concrete analytic realization as locally symmetric
varieties. Even in the early days of modern algebraic geometry the study of
quartic surfaces led to some specific examples of these moduli spaces being
studied in the context of projective geometry. Later advances in complex
analytic and algebraic geometry and in number theory have given us many
very effective tools for studying these varieties and their various
compactifications, and in the last ten years a considerable amount of
progress has been made in understanding the general picture. In this survey
we intend to give a reasonably thorough account of the more recent work,
though mostly without detailed proofs, and to describe sufficiently but not
exhaustively the earlier work of, among others, Satake, Igusa, Mumford and
Tai that has made the recent progress possible.

We confine ourselves to working over the complex numbers. This does not
mean that we can wholly ignore number theory, since much of what is known
depends on interpreting differential forms on Siegel modular varieties as
Siegel modular forms. It does mean, though, that we are neglecting many
important, interesting and difficult questions: in particular, the work of
Faltings and Chai, who extended much of the compactification theory to
$\Spec \ZZ$, will make only a fleeting appearance. To have attempted to
cover this material would have greatly increased the length of this article
and would have led us beyond the areas where we can pretend to competence.

The plan of the article is as follows.

In Section~I we first give a general
description of Siegel modular varieties as complex analytic spaces, and
then explain how to compactify then and obtain projective varieties. There
are essentially two related ways to do this.

In Section~II we start to
understand the birational geometry of these compactified varieties. We
examine the canonical divisor and explain some results which calculate the
Kodaira dimension in many cases and the Chow ring in a few. We also
describe the fundamental group.

In Section~III we restrict ourselves to the
special case of moduli of abelian surfaces (Siegel modular threefolds),
which is of particular interest. We describe a rather general lifting
method, due to Gritsenko in the form we use, which produces Siegel modular
forms of low weight by starting from their behaviour near the boundary of
the moduli space. This enables us to get more precise results about the
Kodaira dimension in a few interesting special cases, due to Gritsenko and
others. Then we describe some results, including a still unpublished
theorem of L.~Borisov, which tend to show that in most cases
the compactified varieties are of general type. In the last part of this
section we examine some finite covers and quotients of moduli spaces of
polarized abelian surfaces, some of which can be interpreted as moduli of
Kummer surfaces. The lifting method gives particularly good results for
these varieties.

In Section~IV we examine three cases, two of them classical, where a Siegel
modular variety (or a near relative) has a particularly good projective
description. These are the Segre cubic and the Burkhardt quartic, which are
classical, and the Nieto quintic, which is on the contrary a surprisingly
recent discovery. There is a huge body of work on the first two and we
cannot do more than summarize enough of the results to enable us to
highlight the similarities among the three cases.

In Section~V we examine the moduli spaces of $(1,t)$-polarized abelian
surfaces (sometimes with level structure) for small~$t$. We begin with the
famous Horrocks-Mumford case, $t=5$, and then move on to the work of
Manolache and Schreyer on $t=7$ and Gross and Popescu on other cases,
especially $t=11$.

In Section~VI we return to the compactification problems and describe very
recent improvements brought about by Alexeev and Nakamura, who (building on
earlier work by Nakamura, Namikawa, Tai and Mumford) have shed some light
on the question of whether there are compactifications of the moduli space
that are really compactifications of moduli, that is, support a proper
universal family.

{\bf Acknowledgements.} Both authors were partially supported by the HCM
network AGE (Algebraic Geometry in Europe), contract no. ERBCHRXCT940557.
They are also grateful to RIMS, Kyoto, for hospitality at different times
during 1996/97 and the first author would like to thank MSRI for
hospitality in the autumn of 1998. We are also grateful to the many people
mentioned in this article who answered our questions about their own work,
and in particular to V.~Alexeev, M.~Gross and S.~Popescu, I.~Nieto, and
N.~Manolache and F.-O.~Schreyer for allowing us access to unpublished
notes.

\section{Siegel modular varieties}\label{I}
In this section we give the basic definitions in connection with
Siegel modular varieties and sketch the construction of the Satake and
toroidal compactifications.

\subsection{Arithmetic quotients of the Siegel upper half plane}\label{I.1}
To any point $\tau$ in the {\em upper half plane}
$$
\HH_1=\{\tau\in \CC\>;\ \  \Im\>\tau>0\}
$$
one can associate a lattice
$$
L_{\tau}=\ZZ\tau+\ZZ
$$
and an {\em elliptic curve}
$$
  E_{\tau}=\CC/L_{\tau}.
$$
Since every elliptic curve arises in this way one obtains a surjective map
$$
  \HH_1\rightarrow \{\mbox {elliptic curves} \}/\mbox{ isomorphism}.
$$
The group $\SL(2,\ZZ)$ acts on $\HH_1$ by
$$\left(\begin{array}{cc}
  a & b\\
  c & d\end{array}\right):\tau\mapsto \frac{a\tau+b}{c \tau +d}
$$
and
$$
  E_{\tau}\cong E_{\tau'}\Leftrightarrow \tau\sim\tau'\mbox{ mod }\SL(2,\ZZ).
$$
Hence there is a bijection
$$
  X^{\circ}(1)=\SL(2,\ZZ)\backslash\HH_1\stackrel{1:1}{\longrightarrow}
  \{\mbox {elliptic curves}
  \}/\mbox{ isomorphism}.
$$
The $j$-function is an $\SL(2,\ZZ)$-invariant function on $\HH_1$ and
defines an isomorphism of Riemann surfaces
$$
  j:X^{\circ}(1)\cong \CC.
$$

An {\em abelian variety} (over the complex numbers $\CC$) is a
$g$-dimensional complex torus $\CC^g/L$ which is a projective variety, i.e.
can be embedded into some projective space $\PP^n$. Whereas every
$1$-dimensional torus $\CC/L$ is an algebraic curve, it is no longer true
that every torus $X=\CC^g/L$ of dimension $g\ge 2$ is projective. This is
the case if and only if $X$ admits a {\em polarization}. There are several
ways to define polarizations. Perhaps the most common definition is that
using Riemann forms. A {\em Riemann form} on $\CC^g$ with respect to the
lattice $L$ is a hermitian form $H\ge 0$ on $\CC^g$ whose imaginary part
$H'=\Im(H)$ is integer-valued on $L$, i.e.  defines an alternating bilinear
form
$$
H':L\otimes L\rightarrow \ZZ.
$$
The $\RR$-linear extension of $H'$ to $\CC^g$ satisfies
$H'(x,y)=H'(ix, iy)$ and determines $H$ by the relation
$$
H(x,y)=H'(ix,y)+iH'(x,y).
$$
$H$ is positive definite if and only if $H'$ is non-degenerate. In this
case $H$ (or equivalently $H'$) is called a {\em polarization}. By the
elementary divisor theorem there exists then a basis of $L$ with respect to
which $H'$ is given by the form
$$
\Lambda=\left(
\begin{array}{cc}
0  &  E\\
-E  &  0\end{array}
\right)\quad, \quad E=\left(
\begin{array}{ccc}
e_1 & & \\
&\ddots &\\
&&e_g
\end{array}
\right),
$$
where the $e_1,\ldots,e_g$ are positive integers such that
$e_1|e_2\ldots |e_g$.  The $g$-tuple $(e_1,\ldots,e_g)$ is uniquely
determined by $H$ and is called the {\em type of the polarization}. If
$e_1=\ldots =e_g=1$ one speaks of a {\em principal polarization}. A ({\em
principally) polarized abelian variety} is a pair $(A, H)$ consisting of a
torus $A$ and a (principal) polarization $H$.

Assume we have chosen a basis of the lattice $L$. If we express each basis
vector of $L$ in terms of the standard basis of $\CC^g$ we obtain a matrix
$\Omega \in M(2g\times g,\CC)$ called a {\em period matrix} of $A$. The fact
that $H$ is hermitian and positive definite is equivalent to
$$
^t\Omega\Lambda^{-1}\Omega = 0, \mbox{ and } i \
^t\Omega\Lambda^{-1}{\bar{\Omega}}>0.
$$
These are the {\em Riemann bilinear relations}. We
consider vectors of $\CC^g$ as row vectors. Using the action of
$\operatorname{GL}(g,\CC)$ on row vectors by right multiplication we can
transform the last $g$ vectors of the chosen basis of $L$ to be $(e_1,
0,\ldots,0),(0,e_2,0,\ldots,0),\ldots,(0,\ldots,0,e_g)$. Then $\Omega$ takes
on the form
$$
\Omega=\Omega_{\tau}=\left(\begin{array}{c}\tau\\E\end{array}\right)
$$
and the Riemann bilinear relations translate into
$$
\tau={^t\tau},\quad \Im \tau>0.
$$
In other words, the complex $(g\times g)$-matrix $\tau$ is an element of
the {\em Siegel space of degree }$g$
$$
\HH_g=\{\tau\in M(g\times g,\CC);\tau={^t\tau},\Im\tau>0\}.
$$
Conversely, given a matrix $\tau\in\HH_g$ we can associate to it the period
matrix $\Omega_{\tau}$ and the lattice $L=L_{\tau}$ spanned by the rows of
$\Omega_{\tau}$. The complex torus $A=\CC^g/L_{\tau}$ carries a Riemann form
given by
$$
H(x,y)=x\Im (\tau)^{-1}\  {^t{\bar y}}.
$$
This defines a polarization of type $(e_1,\ldots, e_g)$. Hence
for every given type of polarization we have a surjection
$$
\HH_g\rightarrow\{(A,H);(A,H)
\mbox{ is an } (e_1,\ldots,e_g)\mbox{-polarized ab.var.}\}/\mbox{ isom.}
$$
To describe the set of these isomorphism classes we have to see what
happens when we change the basis of $L$. Consider the {\em symplectic
group}
$$
\on{Sp}(\Lambda, \ZZ)=\{h \in
\operatorname{GL}(2g,\ZZ);\  h\Lambda^th=\Lambda\}.
$$
As usual we write elements $h\in\on{Sp}(\Lambda,\ZZ)$ in the form
$$
h=\left(
\begin{array}{cc}A & B\\ C & D\end{array}\right);\quad  A,\ldots, D\in
M(g\times
g,\ZZ).
$$
It is useful to work with the ``right projective space $P$ of
$\GL(g,\CC)$'' i.e. the set of all $(2g\times g)$-matrices of rank
$2$ divided out by the equivalence relation
$$
\left(\begin{array}{c}
M_1\\
M_2
\end{array}\right)\sim\left(\begin{array}{c}
M_1 M\\
M_2 M
\end{array}\right) \mbox{  for
any } M\in
\operatorname{GL}(g,\CC).
$$
Clearly $P$ is isomorphic to the Grassmannian
$G=\operatorname{Gr}(g,\CC^{2g}).$
The group
$\on{Sp}(\Lambda,\ZZ)$ acts on $P$ by
$$
\left(
\begin{array}{cc}
A & B\\
C & D \end{array}\right) \left[\begin{array}{c}
M_1\\
M_2
\end{array}\right]=\left[
\begin{array}{c}
AM_1+BM_2\\
CM_1+DM_2
\end{array}\right]
$$
where $[ \ ] $ denotes equivalence classes in $P$. One can embed $\HH_g$ into
$P$ by $\tau\mapsto \left[ \begin{array}{c}\tau\\E\end{array}\right]$. Then
the action of $\on{Sp}(\Lambda, \ZZ)$ restricts to an action on the
image of $\HH_g$ and is given by
$$
\left(
\begin{array}{cc}
A & B\\
C & D \end{array}\right) \left[\begin{array}{c}
\tau\\
E
\end{array}\right]=\left[
\begin{array}{c}
A\tau +BE\\
C\tau +DE
\end{array}\right]=\left[
\begin{array}{c}
(A\tau+BE)(C\tau+DE)^{-1}E\\
E
\end{array}
\right].
$$
In other words, $\on{Sp}(\Lambda,\ZZ)$ acts on $\HH_g$ by
$$
\left(
\begin{array}{cc}
A & B\\
C & D \end{array}\right):\tau\mapsto (A\tau+BE)(C\tau+DE)^{-1}E.
$$
We can then summarize our above discussion with the observation that for a
given type $(e_1,\ldots, e_g)$ of a polarization the quotient
$$
{\cal A}_{e_1,\ldots,e_g}=\on{Sp}(\Lambda, \ZZ)\backslash \HH_g
$$
parametrizes the isomorphism classes of $(e_1,\ldots,e_g)$-polarized
abelian varie\-ties, i.e. ${\cal A}_{e_1,\ldots,e_g}$ is the coarse
moduli space of $(e_1,\ldots,e_g)$-polarized abelian varieties. (Note
that the action of $\on{Sp}(\Lambda,\ZZ)$ on $\HH_g$ depends on the
type of the polarization.)  If we consider principally polarized
abelian varieties, then the form $\Lambda$ is the standard symplectic
form
$$
J=\left(
\begin{array}{cc}
0          &  {\bf 1}_g\\
-{\bf 1}_g &   0
\end{array}
\right)
$$
and $\on{Sp}(\Lambda,\ZZ)=\on{Sp}(2g,\ZZ)$ is the standard symplectic
integer group. In this case we use the notation
$$
{\cal A}_g={\cal A}_{1,\ldots,1}=\on{Sp}(2g,\ZZ)\backslash \HH_g.
$$
This clearly generalizes the situation which we encountered with elliptic
curves. The space $\HH_1$ is just the ordinary upper half plane and
$\on{Sp}(2,\ZZ)=\SL(2,\ZZ)$. We also observe that multiplying the
type of a polarization by a common factor does not change the moduli space.
Instead of the group $\on{Sp}(\Lambda,\ZZ)$ one can also use a
suitable conjugate which is a subgroup of $\on{Sp}(J,\QQ)$. One can then
work with the standard symplectic form and the usual action of the
symplectic group on Siegel space, but the elements of the conjugate
group will in general have rational and no longer just integer entries.

One is often interested in polarized abelian varieties with extra structures,
the so-called {\em level structures}. If $L$ is a lattice equipped with a
non-degenerate form $\Lambda$ the {\em dual lattice} $L^{\vee}$ of $L$ is
defined by
$$
L^{\vee}=\{y\in L\otimes \QQ;\ \Lambda(x,y)\in\ZZ \mbox { for all } x \in L\}.
$$
Then $L^{\vee}/L$ is non-canonically isomorphic to $(\ZZ_{e_1}\times\ldots
\times \ZZ_{e_g})^2$. The group $L^{\vee}/L$ carries a skew form induced by
$\Lambda$ and the group $(\ZZ_{e_1}\times\ldots \times \ZZ_{e_g})^2$ has a
$\QQ/\ZZ$-valued skew form which with respect to the canonical generators is
given by
$$
\left(
\begin{array}{cc}
0       &  E^{-1}\\
-E^{-1} & 0
\end{array}
\right).
$$
If $(A,H)$  is a polarized abelian variety, then a {\em canonical level
structure} on $(A,H)$ is a symplectic isomorphism
$$
\alpha:L^{\vee}/L\rightarrow (\ZZ_{e_1}\times\ldots\times \ZZ_{e_g})^2
$$
where the two groups are equipped with the forms described above. Given
$\Lambda$ we can define the group.
$$
\Sp^\lev(\Lambda,\ZZ):=\{h\in\Sp(\Lambda,\ZZ);\
h|_{L^{\vee}/L}=\mbox{id }_{L^{\vee}/L}\}.
$$
The quotient space
$$
\cA^\lev_{e_1,\ldots,e_g}:=\Sp^\lev(\Lambda,\ZZ)\backslash\HH_g
$$
has the interpretation
$$
\begin{array}{r}
{\cal A}^{\lev}_{e_1,\ldots, e_g}=\{(A,H,\alpha);\ (A,H) \mbox
{ is an } (e_1,\ldots,e_g)\mbox{-polarized abelian}\\
\mbox{variety}, \alpha \mbox { is a canonical level structure}\}/\mbox{ isom}.
\end{array}
$$
If $\Lambda$ is a multiple $n J$ of the standard symplectic form then
$\on{Sp}(n J,\ZZ)=\on{Sp}(J, \ZZ)$ but
$$
\Gamma_g(n):=\on{Sp}^{\lev}(n J,\ZZ)=\{h\in\on{Sp}(J,\ZZ);\
h\equiv{\bf 1}\mbox{ mod } n\}.
$$
This group is called the {\em principal congruence subgroup} of level $n$. A
{\em level}-$n$ {\em structure} on a principally polarized abelian variety
$(A,H)$ is a canonical level structure in the above sense for the polarization
$nH$. The space
$$
{\cal A}_g(n):=\Gamma_g(n)\backslash \HH_g
$$
is the moduli space of principally polarized abelian varieties with a
level-$n$ structure.

The groups $\on{Sp}(\Lambda,\ZZ)$ act properly discontinuously on the
Siegel space $\HH_g$. If $e_1\ge 3$ then
$\Sp^{\lev} (\Lambda,\ZZ)$ acts freely and
consequently the spaces ${\cal
A}^{\lev}_{e_1,\ldots,e_g}$ are smooth in this case. The
finite group
$\Sp(\Lambda,\ZZ)/\Sp^{\lev}(\Lambda,\ZZ)$ acts
on $\cA^{\lev}_{e_1,\ldots, e_g}$ with quotient
${\cal A}_{e_1,\ldots, e_g}$. In particular, these spaces have at most
finite quotient singularities.

A torus $A=\CC^g/L$ is projective if and only if there exists an ample line
bundle ${\cal L}$ on it. By the Lefschetz theorem the first Chern class
defines an isomorphism
$$
c_1:\operatorname{NS}(A)\cong H^2(A,\ZZ)\cap H^{1,1}(A,\CC).
$$
The natural identification $H_1(A,\ZZ)\cong L$ induces isomorphisms
$$
H^2(A,\ZZ)\cong \on{Hom}(\bigwedge\nolimits^2 H_1(A,\ZZ),\ZZ)\cong
\on{Hom}(\bigwedge\nolimits^2 L,\ZZ).
$$
Hence given a line bundle $\cal L$ the first Chern class $c_1({\cal
L})$ can be interpreted as a skew form on the lattice $L$. Let
$H':=-c_1({\cal L})\in\on{Hom}(\bigwedge\nolimits^2 L,\ZZ)$. Since
$c_1({\cal L})$ is a $(1,1)$-form it follows that $H'(x,y)=H'(ix, iy)$
and hence the associated form $H$ is hermitian. The ampleness of
${\cal L}$ is equivalent to positive definiteness of $H$. In this way
an ample line bundle defines, via its first Chern class, a hermitian
form $H$. Reversing this process one can also associate to a Riemann
form an element in $H^2(A,\ZZ)$ which is the first Chern class of an
ample line bundle ${\cal L}$. The line bundle ${\cal L}$ itself is
only defined up to translation. One can also view level structures
from this point of view. Consider an ample line bundle ${\cal L}$
representing a polarization $H$. This defines a map
$$
\begin{array}{cccl}
\lambda: & A & \rightarrow & {\hat A} = \mbox{Pic}^0 A\\
         & x & \mapsto     & t^*_x{\cal L}\otimes {\cal L}^{-1}
\end{array}
$$
where $t_x$ is translation by $x$. The map $\lambda$ depends only on
the polarization, not on the choice of the line bundle ${\cal L}$. If
we write $A=\CC^g/L$ then we have $\Ker\lambda\cong L^{\vee}/L$ and
this defines a skew form on $\Ker\lambda$, the {\em Weil
pairing}. This also shows that $\Ker\lambda$ and the group $(\ZZ_{
e_1} \times\ldots\times \ZZ_{e_g})^2$ are (non-canonically)
isomorphic. We have already equipped the latter group with a skew
form. From this point of view a canonical level structure is then
nothing but a symplectic isomorphism
$$
\alpha:\Ker\lambda\cong (\ZZ_{e_1} \times\ldots \times \ZZ_{e_g})^2.
$$

\subsection{Compactifications of Siegel modular varieties}\label{I.2}

We have already observed that the $j$-function defines an isomorphism of
Riemann
surfaces
$$
j:X^{\circ}(1)=\SL(2,\ZZ)\backslash \HH_1\cong \CC.
$$
Clearly this can be compactified to $X(1)=\PP^1=\CC\cup\{\infty\}.$ It
is, however, important to understand this compactification more
systematically. The action of the group $\SL(2,\ZZ)$ extends to an action
on
$$
{\overline{\HH}}_1=\HH_1\cup\QQ\cup\{i \infty\}.
$$
The extra points $\QQ\cup\{i \infty\}$ form one orbit under this action
and we can set
$$
X(1)=\SL (2,\ZZ)\backslash{\overline{\HH}}_1.
$$
To understand the structure of $X(1)$ as a Riemann surface we have to
consider the stabilizer
$$
P(i \infty)=\left\{
\pm\left(
\begin{array}{cc}
1 & n\\
0 & 1
\end{array}
\right); n\in \ZZ
\right\}
$$
of the point $i \infty$. It acts on $\HH_1$ by $\tau\mapsto \tau+n$.
Taking the quotient by $P(i\infty)$ we obtain the map
$$
\begin{array}{ccl}
\HH_1 & \rightarrow &   D^*_1=\{z\in\CC;\quad 0<|z|<1\}\\
\tau  & \mapsto     &   t=e^{2\pi i\tau}.
\end{array}
$$
Adding the origin gives us the ``partial compactification'' $D_1$ of
$D_1^*$.  For $\varepsilon$ sufficiently small no two points in the
punctured disc $D^*_{\varepsilon}$ of radius $\varepsilon$ are identified
under the map from $D^*_1$ to the quotient $\SL(2,\ZZ)\backslash\HH_1$.
Hence we obtain $X(1)$ by
$$
X(1)=X^{\circ}(1)\cup_{D^*_{\varepsilon}} D_{\varepsilon}.
$$
This process is known as ``adding the cusp $i \infty$''. If we take an
arbitrary arithmetic subgroup $\Gamma\subset \SL(2,\ZZ)$ then
$\QQ\cup\{i\infty\}$ will in general have several, but finitely many,
orbits. However, given a representative of such an orbit we can always
find an element in $\SL(2,\ZZ)$ which maps this representative
to $i \infty$. We can then perform the above construction once more,
the only difference being that we will, in general, have to work
with a subgroup of $P(i\infty)$. Using this process we can always
compactify the quotient $X^{\circ}(\Gamma)=\Gamma\backslash\HH_1$, by
adding a finite number of cusps, to a compact Riemann surface
$X(\Gamma)$.

The situation is considerably more complicated for higher genus $g$
where it is no longer the case that there is a unique compactification
of a quotient ${\cal A}(\Gamma)=\Gamma\backslash\HH_g$. There have
been many attempts to construct suitable compactifications of ${\cal
A}(\Gamma)$. The first solution was given by Satake \cite{Sa} in the
case of ${\cal A}_g$. Satake's compactification ${\bar{\cal A}}_g$ is
in some sense minimal. The boundary ${\bar{\cal A}}_g\backslash{\cal
A}_g$ is set-theoretically the union of the spaces ${\cal A}_i, i\le
g-1$. The projective variety ${\bar{\cal A}}_g$ is normal but highly
singular along the boundary. Satake's compactification was later
generalized by Baily and Borel to arbitrary quotients of symmetric
domains by arithmetic groups. By blowing up along the boundary, Igusa
\cite{I3} constructed a partial desingularization of Satake's
compactification. The boundary of Igusa's compactification has
codimension $1$. The ideas of Igusa together with work of Hirzebruch
on Hilbert modular surfaces were the starting point for Mumford's
general theory of toroidal compactifications of quotients of bounded
symmetric domains \cite{Mu3}. A detailed description of this theory
can be found in \cite{AMRT}. Namikawa showed in \cite{Nam2} that
Igusa's compactification is a toroidal compactification in Mumford's
sense. Toroidal compactifications depend on the choice of cone
decompositions and are, therefore, not unique. The disadvantage of
this is that this makes it difficult to give a good modular
interpretation for these compactifications. Recently, however, Alexeev
and Nakamura \cite{AN},\cite{Ale2}, partly improving work of Nakamura
and Namikawa \cite{Nak1},\cite{Nam1}, have made progress by
showing that the toroidal compactification ${\cal A}_g^*$ which is
given by the second Voronoi decomposition represents a good
functor. We shall return to this topic in chapter VI of our survey article.

This survey article is clearly not the right place to give a complete
exposition of the construction of compactifications of Siegel modular
varieties. Nevertheless we want to sketch the basic ideas behind the
construction of the Satake compactification and of toroidal
compactifications. We shall start with the {\em Satake
compactification}. For this we consider an arithmetic subgroup
$\Gamma$ of $\on{Sp}(2g,\QQ)$ for some $g\ge 2$. (This is no
restriction since the groups $\on{Sp}(\Lambda,\ZZ)$ which arise for
non-principal polarizations are conjugate to subgroups of
$\on{Sp}(2g,\QQ)).$ A {\em modular form} of {\em weight} $k$ with
respect to the group $\Gamma$ is a holomorphic function
$$
F:\HH_g\longrightarrow \CC
$$
with the following transformation behaviour with respect to the group
$\Gamma$:
$$
F(M\tau)=\det(C\tau+D)^k F(\tau)\quad \mbox{for all } M=\left(
\begin{array}{cc}
A & B\\
C & D
\end{array}\right)
\in\Gamma.
$$
(For $g=1$ one has to add the condition that $F$ is holomorpic at the
cusps, but this is automatic for $g\ge 2$). If $\Gamma$ acts freely
then the automorphy factor $\det(C\tau+D)^k$ defines a line bundle
$L^k$ on the quotient $\Gamma\backslash \HH_g$. In general some
elements in $\Gamma$ will have fixed points, but every such element is
torsion and the order of all torsion elements in $\Gamma$ is bounded
(see e.g. \cite[p.120]{LB}). Hence, even if $\Gamma$ does not act
freely, the modular forms of weight $n k_0$ for some suitable integer
$k_0$ and $n\ge 1$ are sections of a line bundle $L^{nk_0}$. The space
$M_k(\Gamma)$ of modular forms of fixed weight~$k$ with respect to
$\Gamma$ is a finite-dimensional vector space and the elements of
$M_{nk_0}(\Gamma)$ define a rational map to some projective space
$\PP^N$. If $n$ is sufficiently large it turns out that this map is
actually an immersion and the Satake compactification ${\cal
A}(\Gamma)$ can be defined as the projective closure of the image of
this map.

There is another way of describing the Satake compactification which
also leads us to toroidal compactifications. The {\em Cayley
transformation}
$$
\begin{array}{rcl}
\Phi:\HH_g\  &  \rightarrow  &  \on{Sym} (g,\CC)\\
\tau           &  \mapsto      &  (\tau-i{\bf 1})(\tau+i{\bf 1)}^{-1}
\end{array}
$$
realizes $\HH_g$ as the symmetric domain
$$
{\cal D}_g=\{Z\in\on{Sym}(g,\CC);\  {\bf 1}-Z{\bar Z}>0\}.
$$
Let ${\bar{\cal D}}_g$ be the topological closure of ${\cal D}_g$ in
$\Sym(g,\CC)$. The action of $\on{Sp}(2g,\RR)$ on $\HH_g$
defines, via the Cayley transformation, an action on ${\cal D}_g$
which extends to ${\bar{\cal D}}_g$. Two points in ${\bar{\cal D}}_g$
are called equivalent if they can be connected by finitely many
holomorphic curves. Under this equivalence relation all points in
${\cal D}_g$ are equivalent. The equivalence classes of ${\bar{\cal
D}}_g\backslash {\cal D}_g$ are called the proper boundary components
of ${\bar{\cal D}}_g$. Given any point $Z\in {\bar{\cal D}}_g$ one can
associate to it the real subspace $U(Z)=\ker\psi(Z)$ of $\RR^{2g}$
where
$$
\psi(Z):\RR^{2g}\rightarrow \CC^g,\nu\mapsto \nu\left(
\begin{array}{c}
i({\bf 1} +Z)\\
{\bf 1}-Z
\end{array}
\right)
$$
Then $U(Z)$ is an isotropic subspace of $\RR^{2g}$ equipped with the
standard symplectic form $J$. Moreover $U(Z)\neq 0$ if and only if
$Z\in{\bar{\cal D}}_g\backslash {\cal D}_g$ and $U(Z_1)=U(Z_2)$ if and
only if $Z_1$ and $Z_2$ are equivalent. This defines a bijection
between the proper boundary components of ${\bar{\cal D}}_g$ and the
non-trivial isotropic subspaces of $\RR^{2g}$. A boundary component
$F$ is called {\em rational} if its stabilizer subgroup $P(F)$ in
$\on{Sp}(2g,\RR)$ is defined over the rationals or, equivalently, if
$U(F)$ is a rational subspace, i.e. can be generated by rational
vectors. Adding the rational boundary components to ${\cal D}_g$ one
obtains the {\em rational closure} ${\cal D}_g^{\rm rat}$ of ${\cal
D}_g$. This can be equipped with either the Satake topology or the
cylindrical topology. The Satake compactification, as a topological
space, is then the quotient $\Gamma\backslash{\cal D}_g^{\rm
rat}$. (The Satake topology and the cylindrical topology are actually
different, but the quotients turn out to be homeomorphic.) For $g=1$
the above procedure is easily understood: the Cayley transformation
$\psi$ maps the upper half plane $\HH_1$ to the unit disc $D_1$.
Under this transformation the rational boundary points $\QQ\cup\{i
\infty\}$ of $\HH_1$ are mapped to the rational boundary points of
$D_1$. The relevant topology is the image under $\psi$ of the {\em
horocyclic topology} on $\overline{\HH}_1=\HH_1\cup\QQ\cup\{i
\infty\}$.

Given two boundary components $F$ and $F'$ with $F\neq F'$ we say that
$F$ is {\em adjacent} to $F'$ (denoted by $F'\succ F$) if $F\subset
{\overline{ F'}}$. This is the case if and only if $U(F')\subsetneqq
U(F)$. In this way we obtain two partially ordered sets, namely
$$
\begin{array}{ccl}
(X_1,<) & =  & (\{\mbox{proper rational boundary components } F \mbox{ of }
{\cal D}_g
\}, \succ)\\
(X_2,<) & =  & (\{\mbox{non-trivial isotropic subspaces }U\mbox{ of }\QQ^g
\}, \subsetneqq).
\end{array}
$$
The group $\on{Sp}(2g,\QQ)$ acts on both partially ordered sets as a
group of automorphisms and the map $f:X_1\rightarrow X_2$ which
associates to each $F$ the isotropic subspace $U(F)$ is an
$\on{Sp}(2g,\QQ)$-equivariant isomorphism of partially ordered
sets. To every partially ordered set $(X,<)$ one can associate its
{\em simplicial realization} $\on{SR}(X)$ which is the simplicial
complex consisting of all simplices $(x_0,\ldots, x_n)$ where
$x_0,\ldots,x_n \in X$ and $ x_0<x_1<\ldots<x_n$. The {\em Tits
building} ${\cal T}$ of $\on{Sp}(2g,\QQ)$ is the simplicial complex
${\cal T}=\on{SR}(X_1)=\on{SR}(X_2)$. If $\Gamma$ is an arithmetic
subgroup of $\on{Sp}(2g,\QQ)$, then the Tits building of $\Gamma$ is
the quotient ${\cal T}(\Gamma)=\Gamma\backslash{\cal T}$.

For any boundary component $F$ we can define its stabilizer in
$\Sp(2g,\RR)$ by
$$
{\cal P}(F)=\{h\in \on{Sp}(2g,\RR);\  h(F)=F\}.
$$
If $U=U(F)$ is the associated isotropic subspace, then
$$
{\cal P}(F)={\cal P}(U)=\{h\in \on{Sp}(2g,\RR);\  Uh^{-1} =U\}.
$$
The group ${\cal P}(F)$ is a maximal parabolic subgroup of
$\Sp(2g,\RR)$. More generally, given any flag
$U_1\subsetneqq\ldots\subsetneqq U_l$ of isotropic subspaces, its
stabilizer is a parabolic subgroup of $\on{Sp}(2g,\RR)$. Conversely
any parabolic subgroup is the stabilizer of some isotropic flag. The
maximal length of an isotropic flag in $\RR^{2g}$ is $g$ and the
corresponding subgroups are the minimal parabolic subgroups or {\em
Borel subgroups} of $\on{Sp}(2g,\RR)$. We have already remarked that a
boundary component $F$ is rational if and only if the stabilizer
${\cal P}(F)$ is defined over the rationals, which happens if and only
if $U(F)$ is a rational subspace. More generally an isotropic flag is
rational if and only if its stabilizer is defined over $\QQ$. This
explains how the Tits building ${\cal T}$ of $\on{Sp}(2g,\QQ)$ can be
defined using parabolic subgroups of $\on{Sp}(2g,\RR)$ which are
defined over $\QQ$. The Tits building of an arithmetic subgroup
$\Gamma$ of $\on{Sp}(2g,\QQ)$ can, therefore, also be defined in terms
of conjugacy classes of groups $\Gamma\cap{\cal P}(F)$.

As an example we consider the integer symplectic group $\on{Sp}(2g,\ZZ)$.
There exists exactly one maximal isotropic flag modulo the action of
$\on{Sp}(2g,\ZZ)$, namely
$$
\{0\}\subsetneqq U_1\subsetneqq U_2\subsetneqq\ldots\subsetneqq U_g;\quad
U_i=\on{span }(e_1,\ldots,e_i).
$$
Hence the Tits building ${\cal T}(\on{Sp}(2g,\ZZ))$ is a
$(g-1)$-simplex whose vertices correspond to the space $U_i$. This
corresponds to the fact that set-theoretically
$$
{\bar{\cal A}}_g={\cal A}_g\ {\amalg}\ {\cal A}_{g-1}\ {\amalg}\ \ldots
\ {\amalg}\ {\cal A}_1\ {\amalg}\ {\cal A}_0.
$$

With these preparations we can now sketch the construction of a
toroidal compactification of a quotient ${\cal
A}(\Gamma)=\Gamma\backslash\HH_g$ where $\Gamma$ is an arithmetic
subgroup of $\on{Sp}(2g,\QQ)$. We have to compactify ${\cal
A}(\Gamma)$ in the direction of the cusps, which are in 1-to-1
correspondence with the vertices of the Tits building ${\cal
T}(\Gamma)$. We shall first fix one cusp and consider the associated
boundary component $F$, resp. the isotropic subspace $U=U(F)$. Let
${\cal P}(F)$ be the stabilizer of $F$ in $\on{Sp}(2g,\RR)$. Then
there is an exact sequence of Lie groups
$$
1\rightarrow {\cal P}'(F)\rightarrow {\cal P}(F)\rightarrow{\cal
P}''(F)\rightarrow 1
$$
where ${\cal P}'(F)$ is the centre of the unipotent radical $R_u({\cal
P}(F))$ of ${\cal P}(F)$. Here ${\cal P}'(F)$ is a real vector space
isomorphic to $\Sym(g',\RR)$ where $g'=\dim U(F)$. Let
$P(F)={\cal P}(F)\cap\Gamma, P'(F)={\cal P}'(F)\cap\Gamma$ and
$P''(F)=P(F)/P'(F)$. The group $P'(F)$ is a lattice of maximal rank in
${\cal P}'(F)$. To $F$ one can now associate a torus bundle ${\cal
X}(F)$ with fibre $T=P'(F)\otimes_{\ZZ}\CC/P'(F)\cong(\CC^*)^{g'}$
over the base $S=F\times V(F)$ where $V(F)=R_u({\cal P}(F))/{\cal
P}'(F)$ is an abelian Lie group and hence a vector space. To construct
a partial compactification of ${\cal A}(\Gamma)$ in the direction of
the cusp corresponding to $F$, one then proceeds as follows:
\begin{enumerate}
\item[(1)] Consider the partial quotient $X(F)=P'(F)\backslash\HH_g$. This
  is a torus bundle with fibre $(\CC^*)^{g'}$ over some open subset of
  $\CC^{\frac 12 g(g+1)-g'}$ and can be regarded as an open subset of the
  torus bundle ${\cal X}(F)$.
\item[(2)] Choose a fan $\Sigma$ in the real vector space ${\cal
  P}'(F)\cong \on{Sym}(g',\RR)$ and construct a trivial bundle ${\cal
  X}_{\Sigma}(F)$ whose fibres are torus embeddings.
\item[(3)] If $\Sigma$ is chosen compatible with the action of $P''(F)$,
  then the action of $P''(F)$ on ${\cal X}(F)$ extends to an action of
  $P''(F)$ on ${\cal X}_{\Sigma}(F)$.
\item[(4)] Denote by $X_{\Sigma}(F)$ the interior of the closure of $X(F)$
  in ${\cal X}_{\Sigma}(F)$. Define the partial compactification of ${\cal
  A}(\Gamma)$ in the direction of $F$ as the quotient space
  $Y_{\Sigma}(F)=P''(F)\backslash X_{\Sigma}(F)$.
\end{enumerate}
To be able to carry out this programme we may not choose the fan $\Sigma$
arbitrarily, but we must restrict ourselves to {\em admissible} fans
$\Sigma$ (for a precise definition see [Nam2, Definition 7.3]). In
particular $\Sigma$ must define a cone decomposition of the cone
$\on{Sym}_+(g',\RR)$ of positive definite symmetric $(g'\times
g')$-matrices. The space $Y_{\Sigma}(F)$ is called the partial
compactification in the direction $F$.

The above procedure describes how to compactify ${\cal A}(\Gamma)$ in the
direction of one cusp $F$. This programme then has to be carried out for
each cusp in such a way that the partial compactifications glue together
and give the desired toroidal compactification. For this purpose we have to
consider a collection ${\tilde{\Sigma}}=\{\Sigma(F)\}$ of fans
$\Sigma(F)\subset{\cal P}'(F)$. Such a collection is called an {\em
admissible collection of fans} if
\begin{enumerate}
\item[(1)]
Every fan $\Sigma(F)\subset{\cal P}'(F)$ is an admissible fan.
\item[(2)]
If $F=g(F')$ for some $g\in\Gamma$, then $\Sigma(F)=g(\Sigma(F'))$ as fans in
the space ${\cal P}'(F)=g({\cal P}'(F'))$.
\item[(3)]
If $F'\succ F$ is a pair of adjacent rational boundary components, then
equality
$\Sigma(F')=\Sigma(F)\cap{\cal P}'(F')$ holds as fans in ${\cal
P}'(F')\subset{\cal P}'(F)$.
\end{enumerate}
The conditions (2) and (3) ensure that the compactifications in the
direction of the various cusps are compatible and can be glued together.
More precisely we obtain the following:

\begin{enumerate}
\item[($2'$)] If $g\in\Gamma$ with $F=g(F')$, then there exists a natural
  isomorphism ${\tilde g}:X_{\Sigma(F')}(F')\rightarrow X_{\Sigma(F)}(F)$.
\item[($3'$)] Suppose $F'\succ F$ is a pair of adjacent rational boundary
  components. Then $P'(F')\subset P'(F)$ and there exists a natural
  quotient map $\pi_{0}(F', F):X(F')\rightarrow X(F)$. Because of (3) this
  extends to an \'etale map: $\pi(F',F):X_{\Sigma(F')}(F')\rightarrow
  X_{\Sigma(F)}(F)$.
\end{enumerate}
We can now consider the disjoint union
$$
X=\coprod\limits_F X_{\Sigma(F)} (F)
$$
over all rational boundary components~$F$. One can define an
equivalence relation on $X$ as follows: if $x\in X_{\Sigma(F)}(F)$ and
$x'\in X_{\Sigma(F')}(F')$, then
\begin{enumerate}
\item[(a)]
$x\sim x'$ if there exists $g\in \Gamma$ such that $F=g(F')$ and $x={\tilde
g}(x')$.
\item[(b)]
$x\sim x'$ if $F'\succ F$ and $\pi(F', F)(x')=x$.
\end{enumerate}

The {\em toroidal compactification} of ${\cal A}(\Gamma)$ defined by the
admissible collection of fans ${\tilde{\Sigma}}$ is then the space
$$
{\cal A}(\Gamma)^{\ast}=X/\sim.
$$
Clearly ${\cal A}(\Gamma)^{\ast}$ depends on ${\tilde\Sigma}$. We
could also have described ${\cal A}(\Gamma)^{\ast}$ as $Y/\sim$ where
$Y=\amalg\ Y_{\Sigma(F)}(F)$ and the equivalence relation $\sim$ on
$Y$ is induced from that on $X$. There is a notion of a {\em
projective} admissible collection of fans (see [Nam2, Definition
7.22]) which ensures that the space ${\cal A}(\Gamma)^{\ast}$ is
projective.

For every toroidal compactification there is a natural map $\pi:{\cal
A}(\Gamma)^{\ast}\rightarrow {\bar{{\cal A}}}(\Gamma)$ to the
Satake compactification. Tai, in \cite{AMRT}, showed that if ${\cal
A}(\Gamma)^{\ast}$ is defined by a projective admissible collection of
fans, then $\pi$ is the normalization of the blow-up of some ideal
sheaf supported on the boundary of ${\bar{{\cal A}}}(\Gamma)$.

There are several well known cone decompositions for
$\on{Sym}_+(g',\RR)$: see e.g. [Nam2, section 8]. The {\em central
cone decomposition} was used by Igusa \cite{I1} and leads to the {\em
Igusa compactification}. The most important decomposition for our
purposes is the {\em second Voronoi decomposition}. The corresponding
compactification is simply called the {\em Voronoi
compactification}. The Voronoi compactification ${\cal
A}(\Gamma)^{\ast}={\cal A}_g^{\ast}$ for $\Gamma=\on{Sp}(2g,\ZZ)$ is a
projective variety \cite{Ale2}. For $g=2$ all standard known cone
decompositions coincide with the {\em Legendre decomposition}.

\section{Classification theory}\label{II}

Here we discuss known results about the Kodaira dimension of Siegel
modular varieties and about canonical and minimal models. We also
report on some work on the fundamental group of Siegel modular
varieties.

\subsection{The canonical divisor}\label{II.1}
If one wants to prove results about the Kodaira dimension of Siegel modular
varieties, one first has to understand the canonical divisor. For an
element $\tau \in \HH_g$ we write
$$
\tau = \left(
\begin{array}{ccc|c}
\tau_{11}    & \cdots & \tau_{1, g-1} & \tau_{1g}\\
\vdots       &        &               & \vdots\\
\tau_{1,g-1} &\cdots  & \tau_{g-1,g-1}& \tau_{g-1, g}\\
\hline
\tau_{1g}       &\cdots  & \tau_{g-1, g} & \tau_{gg}
\end{array}
\right)=\left(
\begin{array}{ccc|c}
{}&{}&{}&{}\\
{}&\tau'&{}& ^t z\\
{}&{}&{}&{}\\
\hline
{}&z{}& & \tau_{gg}
\end{array}
\right).
$$
Let
$$
d\tau=d\tau_{11} \wedge d\tau_{12}\wedge\ldots\wedge d\tau_{gg}.
$$
If $F$ is a modular form of weight $g+1$ with respect to an arithmetic
group $\Gamma$, then it is easy to check that the form $\omega=Fd\tau$
is $\Gamma$-invariant. Hence, if $\Gamma$ acts {\em freely}, then
$$
K_{{\cal A}(\Gamma)} = (g+1)L
$$
where $L$ is the line bundle of modular forms, i.e. the line bundle
given by the automorphy factor $\det(C\tau+D)$. If $\Gamma$ does not act
freely, let $\opn{\cal A} (\Gamma)={\cal A}(\Gamma)\backslash R$ where $R$
is the branch locus of the quotient map $\HH_g\rightarrow {\cal
A}(\Gamma)$.  Then by the above reasoning it is still true that
$$
K_{\opn{\cal A}(\Gamma)} = (g+1)L|_{\opn{\cal A}(\Gamma)}.
$$
In order to describe the canonical bundle on a toroidal compactification
${\cal A}(\Gamma)^{\ast}$ we have to understand the behaviour of the
differential form $\omega$ at the boundary. To simplify the exposition, we
shall first consider the case $\Gamma_g=\on{Sp}(2g,\ZZ)$. Then there
exists, up to the action of $\Gamma$, exactly one maximal boundary
component $F$. We can assume that $U(F)=U=\on{span}(e_g)$. The stabilizer
$P(F)=P(U)$ of $U$ in $\Gamma_g$ is generated by elements of the form
$$
\begin{array}{lr}
{g_1=\left(
\begin{array}{cccc}
A & 0 & B & 0\\
0 & 1 & 0 & 0\\
C & 0 & D & 0\\
0 & 0 & 0 & 1
\end{array}
\right),}
&
{g_2=\left(
\begin{array}{cccc}
{\bf 1}_{g-1} & 0 & 0 & 0\\
0 & \pm1 & 0 & 0\\
0 & 0 & {\bf 1}_{g-1} & 0\\
0 & 0 & 0 & \pm1
\end{array}
\right),}
\end{array}
$$
$$
\begin{array}{lr}
{g_3=\left(
\begin{array}{cccc}
{\bf 1}_{g-1} & 0 & 0 & ^tN\\
M & 1 & N & 0\\
0 & 0 & {\bf 1}_{g-1} & -^tM\\
0 & 0 & 0 & 1
\end{array}
\right),}
&
{g_4=\left(
\begin{array}{cccc}
{\bf 1}_{g-1} & 0 & 0 & 0\\
0 & 1 & 0 & S\\
0 & 0 & {\bf 1}_{g-1} & 0\\
0 & 0 & 0 & 1
\end{array}
\right),}
\end{array}
$$
where $
\left(
\begin{array}{cc}
A & B\\
C & D
\end{array}
\right) \in \Gamma_{g-1}$, $M,N\in \ZZ^{g-1}$ and $S\in \ZZ$.

The group $P'(F)$ is the rank 1 lattice generated by $g_4$, and the
partial quotient with respect to $P'(F)$ is given by
$$
\begin{array}{rcl}
e(F):\ \HH_g & \longrightarrow & \HH_{g-1} \times \CC^{g-1}\times
\CC^{\ast} \\
\tau & \longmapsto & (\tau', z, t=e^{2\pi i\tau_{gg}}).
\end{array}
$$
Here $\HH_{g-1}\times \CC^{g-1}\times \CC^{\ast}$ is a rank 1 torus
bundle over $\HH_{g-1}\times\CC^{g-1}=F\times V(F)$. Partial
compactification in the direction of $F$ consists of adding
$\HH_{g-1}\times \CC^{g-1}\times\{0\}$ and then taking the quotient
with respect to $P''(F)$. Since $d\tau_{gg}=(2\pi i)^{-1} dt/t$ it
follows that
$$
\omega=(2\pi i)^{-1} F\frac{d\tau_{11}\wedge\ldots\wedge d\tau_{g-1,
g}\wedge dt}{t}
$$
has a pole of order $1$ along the boundary, unless $F$ vanishes there.
Moreover, since $F(g_4(\tau))=F(\tau)$ it follows that $F$ has a
Fourier expansion
$$
F(\tau)=\sum\limits_{n\ge 0} F_n(\tau',z)t^n.
$$
A modular form $F$ is a {\em cusp form} if $F_0(\tau',z)=0$, i.e. if
$F$ vanishes along the boundary. (If $\Gamma$ is an arbitrary
arithmetic subgroup of $\on{Sp}(2g,\QQ)$ we have in general several
boundary components and then we require vanishing of $F$ along each of
these boundary components.) The above discussion can be interpreted as
follows. First assume that $\Gamma$ is neat (i.e. the subgroup of
$\CC^*$ generated by the eigenvalues of all elements of $\Gamma$ is
torsion free) and that ${\cal A}(\Gamma)^{\ast}$ is a smooth
compactification with the following property: for every point in the
boundary there exists a representative $x\in X_{\Sigma(F)}(F)$ for
some boundary component such that $X_{\Sigma(F)}(F)$ is smooth at $x$
and $P''(F)$ acts freely at $x$. (Such a toroidal compactification
always exists if $\Gamma$ is neat.) Let $D$ be the boundary divisor of
${\cal A}(\Gamma)^*$. Then
$$
K_{{\cal A}(\Gamma)^*}=(g+1)L-D.
$$
Here $L$ is the extension of the line bundle on modular forms on
${\cal A}(\Gamma)$ to ${\cal A}(\Gamma)^*$. This makes sense since by
construction the line bundle extends to the Satake compactification
${\bar{{\cal A}}}(\Gamma)$ and since there is a natural map
{${\pi:{\cal A}(\Gamma)^{\ast}\rightarrow{\bar{{\cal
A}}}(\Gamma)}$}. We use the same notation for $L$ and
$\pi^{\ast}L$. If $\Gamma$ does not act freely we can define the open
set ${\opn{\cal A}}(\Gamma)^{\ast}$ consisting of
$\opn{\cal A}(\Gamma)$ and those points in the boundary which
have a representative $x\in X_{\Sigma(F)}(F)$ where $P''(F)$ acts
freely at $x$. In this case we still have
$$
K_{\opn{\cal A}(\Gamma)^{\ast}} = ((g+1)L-D)|_{\opn{\cal
A}(\Gamma)^{\ast}}.
$$
This shows in particular that every cusp form $F$ of weight $g+1$ with
respect to $\Gamma$ defines via $\omega=Fd\tau$ a differential
$N$-form on ${\opn{\cal A}}(\Gamma)^{\ast}$ where
$N=\frac{g(g+1)}{2}$ is the dimension of ${\cal A}(\Gamma)$. It is a
non-trivial result of Freitag that every such form can be extended to
any smooth projective model of ${\cal A}(\Gamma)$. If we denote by
$S_k(\Gamma)$ the space of cusp forms of weight $k$ with respect to
$\Gamma$, then we can formulate Freitag's result as follows.

\begin{theorem}[\cite{F}]\label{theo1}
  Let $\tilde{\cal A}(\Gamma)$ be a smooth projective model of ${\cal
  A}(\Gamma)$. Then every cusp form $F$ of weight $g+1$ with respect to
  $\Gamma$ defines a differential form $\omega=Fd\tau$ which extends to
  $\tilde \calA(\Gamma)$. In particular, there is a natural isomorphism
$$
\Gamma(\tilde{\cal A}(\Gamma),
\omega_{\tilde{\cal A}(\Gamma)})\cong S_{g+1}(\Gamma)
$$
and hence $p_g(\tilde{\cal A}(\Gamma))=\dim S_{g+1}(\Gamma)$.
\end{theorem}
\begin{Proof}
See \cite[Satz III.2.6]{F} and the remark following this.
\end{Proof}
Similarly a form of weight $k(g+1)$ which vanishes of order $k$ along the
boundary defines a $k$-fold differential form on $\opn{\cal
A}(\Gamma)^{\ast}$. In general, however, such a form does not extend to a
smooth model $\tilde{\cal A}(\Gamma)$ of ${\cal A}(\Gamma)$.

\subsection{The Kodaira dimension of ${\cal A}_g(n)$}\label{II.2}
\setcounter{theorem}{0}
By the {\em Kodaira dimension} of a Siegel modular variety ${\cal
A}(\Gamma)$ we mean the Kodaira dimension of a smooth projective model of
${\cal A}(\Gamma)$. Such a model always exists and the Kodaira dimension is
independent of the specific model chosen. It is a well known result that
${\cal A}_g$ is of general type for $g\ge 7$. This was first proved by Tai
for $g\ge 9$ \cite{T1} and then improved to $g\ge 8$ by Freitag \cite{F}
and to $g\ge 7$ by Mumford \cite{Mu4}. In this section we want to discuss
the proof of the following result.
\begin{theorem}[\cite{T1},\cite{F},\cite{Mu4},\cite{H2}]\label{theo2}
${\cal A}_g(n)$ is of general type for the following values of $g$ and $n\ge
n_0$:
$$
\begin{array}{c|cccccr}
g   & 2 & 3 & 4 & 5 & 6 & \ge 7\\
\hline
n_0 & 4 & 3 & 2 & 2 & 2 & 1
\end{array}.
$$
\end{theorem}
We have already seen that the construction of differential forms is closely
related to the existence of cusp forms. Using Mumford's extension of
Hirzebruch proportionality to the non-compact case and the Atiyah-Bott
fixed point theorem it is not difficult to show that the dimension of the
space of cusp forms of weight $k$ grows as follows:
$$
\dim S_k (\Gamma_g)\sim 2^{-N-g}k^N V_g {\pi^{-N}}
$$
where
$$
N=\frac{g(g+1)}{2}=\dim {\cal A}_g(n)
$$
and $V_g$ is Siegel's symplectic volume
$$
V_g=2^{g^2+1}\pi^N \prod\limits^g_{j=1}\frac{(j-1)!}{2j!}B_j.
$$
Here $B_j$ are the Bernoulli numbers.

Every form of weight $k(g+1)$ gives rise to a $k$-fold differential form on
$\opn{\cal A}_g(n)$. If $k=1$, we have already seen that these forms extend
by Freitag's extension theorem to every smooth model of ${\cal A}_g(n)$.
This is no longer automatically the case if $k\ge 2$. Then one encounters
two types of obstructions: one is extension to the boundary (since we need
higher vanishing order along $D$), the other type of obstruction comes from
the singularities, or more precisely from those points where $\Gamma_g(n)$
does not act freely. These can be points on ${\cal A}_g(n)$ or on the
boundary. If $n\ge 3$, then $\Gamma_g(n)$ is neat and in particular it acts
freely. Moreover we can choose a suitable cone decomposition such that the
corresponding toroidal compactification is smooth. In this case there are
no obstructions from points where $\Gamma_g(n)$ does not act freely. If
$n=1$ or $2$ we shall, however, always have such points. It is one of the
main results of Tai \cite[Section 5]{T1} that for $g\ge 5$ all resulting
singularities are {\em canonical}, i.e. give no obstructions to extending
$k$-fold differential forms to a smooth model. The remainder of the proof
of Tai then consists of a careful analysis of the obstructions to the
extension of $k$-forms to the boundary. These obstructions lie in a vector
space which can be interpreted as a space of Jacobi forms on
$\HH_{g-1}\times\CC^{g-1}$. Tai gives an estimate of this space in
\cite[Section 2]{T1} and compares it with the dimension formula for
$S_k(\Gamma_g)$.

The approach developed by Mumford in \cite{Mu4} is more geometric in nature.
First recall that
\begin{equation}
K|_{\opn{\cal A}_g^{\ast}(n)} = ((g+1) L-D)|_{\opn{\cal
A}_g^{\ast}(n)}.
\end{equation}
Let $\bar{\Theta}_{\mbox{\scriptsize{null}}}$ be the closure of the
locus of pairs $(A,\Theta)$ where $A$ is an abelian variety and
$\Theta$ is a symmetric divisor representing a principal polarization
such that $\Theta$ has a singularity at a point of order $2$. Then
one can show that for the class of $\bar{\Theta}_{\mbox{\scriptsize
null}}$ on ${\cal A}_g^{\ast}(n)$:
\begin{equation}
[\bar{\Theta}_{\mbox{\scriptsize null}}]=2^{g-2}(2^g+1) L-2^{2g-5}D.
\end{equation}
One can now use (2) to eliminate the boundary $D$ in (1). Since the natural
quotient ${\cal A}_g^{\ast}(n)\rightarrow {\cal A}_g^{\ast}$ is branched of
order $n$ along $D$ one finds the following formula for $K$:
\begin{equation}
K|_{\opn{\cal A}_g^{\ast}(n)}=\left( (g+1)-\frac{2^{g-2}(2^g+1)}{n
2^{2g-5}}\right) L+\frac 1{n 2^{2g-5}}[\bar{\Theta}_{\mbox{\scriptsize null}}].
\end{equation}
In view of Tai's result on the singularities of ${\cal A}^{\ast}_g(n)$
this gives general type whenever the factor in front of $L$ is
positive and $n\ge 3$ or $g\ge 5$. This gives all cases in the list
with two exceptions, namely $(g,n)=(4,2)$ and $(7,1)$. In the first
case the factor in front of $L$ is still positive, but one cannot
immediately invoke Tai's result on canonical singularities. As
Salvetti Manni has pointed out, one can, however, argue as follows. An
easy calculation shows that for every element $\sigma\in\Gamma_g(2)$
the square $\sigma^2\in\Gamma_g(4)$. Hence if $\sigma$ has a fixed
point then $\sigma^2=1$ since $\Gamma_g(4)$ acts freely. But now one
can again use Tai's extension result (see \cite[Remark after Lemma
4.5]{T1} and \cite[Remark after Lemma 5.2]{T1}).

This leaves the case $(g,n)=(7,1)$ which is the main result of \cite{Mu4}.
Mumford considers the locus
$$
N_0=\{(A,\Theta)\>;\ \on{Sing}\Theta\neq\emptyset\}
$$
in ${\cal A}_g$. Clearly this contains $\Theta_{\mbox{\scriptsize null}}$, but
is bigger than $\Theta_{\mbox{\scriptsize null}}$ if $g\ge 4$. Mumford
shows that the class of the closure $\bar{N}_0$ on ${\cal A}^{\ast}_g$ is
\begin{equation}
[\bar{N}_0]=\left(\frac{(g+1)!}{2}+g!\right)L-\frac{(g+1)!}{12}D
\end{equation}
and hence one finds for the canonical divisor:
$$
K|_{\opn{\cal
A}^{\ast}_g(n)}=\frac{12(g^2-4g-17)}{g+1}L+\frac{12}{(g+1)!}[\bar{N}_0].
$$
Since the factor in front of $L$ is positive for $g=7$ one can once more use
Tai's extension result to prove the theorem for $(g,n)=(7,1)$.

The classification of the varieties ${\cal A}_g(n)$ with respect to
the Kodaira dimension is therefore now complete with the exception of
one important case:
\begin{problem}
Determine the Kodaira dimension of ${\cal A}_6$.
\end{problem}
All other varieties ${\cal A}_g(n)$ which do not appear in the above
list are known to be either rational or unirational. Unirationality of
${\cal A}_5$ was proved by Donagi \cite{D} and independently by Mori
and Mukai \cite{MM} and Verra \cite{V}. Unirationality of ${\cal A}_4$
was shown by Clemens \cite{Cl} and unirationality of ${\cal A}_g, g\le
3$ is easy. For $g=3$ there exists a dominant map from the space of
plane quartics to ${\cal M}_3$ which in turn is birational to ${\cal
A}_3$. For $g=2$ one can use the fact that ${\cal M}_2$ is birational
to ${\cal A}_2$ and that every genus $2$ curve is a 2:1 cover of
$\PP^1$ branched in 6 points. Rationality of these spaces is a more
difficult question. Igusa \cite{I1} showed that ${\cal A}_2$ is
rational. The rationality of ${\cal M}_3$, and hence also of ${\cal
A}_3$, was proved by Katsylo \cite{K}. The space ${\cal A}_3(2)$ is
rational by the work of van Geemen \cite{vG} and Dolgachev and Ortland
\cite{DO}. The variety ${\cal A}_2(3)$ is birational to the Burkhardt
quartic in $\PP^4$ and hence also rational. This was proved by Todd in
1936 \cite{To} and Baker in 1942 (see \cite{Ba2}), but see also the
thesis of Finkelnberg \cite{Fi}. The variety ${\cal A}_2(2)$ is
birational to the Segre cubic (cf. \cite{vdG1}) in $\PP^4$ and hence
also rational. The latter two cases are examples of Siegel modular
varieties which have very interesting projective models. We will come
back to this more systematically in chapter IV. It should also be
noted that Yamazaki \cite{Ya} was the first to prove that ${\cal
A}_2(n)$ is of general type for $n\ge 4$.

All the results discussed above concern the case of principal
polarization. The case of non-principal polarizations of type
$(e_1,\ldots, e_g)$ was also studied by Tai.

\begin{theorem}[\cite{T2}]\label{theo3}
The moduli space ${\cal A}_{e_1,\ldots,e_g}$ of abelian varieties with a
polarization of type $(e_1,\ldots,e_g)$ is of general type if
either $g\ge 16$ or $g\ge 8$ and all $e_i$ are odd and sums of two squares.
\end{theorem}

The essential point in the proof is the construction of sufficiently many cusp
forms with  high vanishing order along the boundary. These modular forms are
obtained as pullbacks of theta series on Hermitian or quaternionic upper half
spaces.

More detailed results are known in the case of abelian surfaces $(g=2)$. We
will discuss this separately in chapters III and~V.

By a different method, namely using symmetrization of modular forms, Gritsenko
has shown the following:

\begin{theorem}[\cite{Gr1}]\label{theo4}
For every integer $t$ there is an integer $g(t)$ such that the moduli space
$\calA_{1,\ldots,1,t}$ is of general type for $g\geq g(t)$. In particular
$\calA_{1,\ldots,1,2}$ is of general type for $g\geq 13$.
\end{theorem}

\begin{Proof}
See [Gr1, Satz 1.1.10], where an explicit bound for $g(t)$ is given.
\end{Proof}

Once one has determined that a variety is of general type it is
natural to ask for a minimal or canonical model. For a given model
this means asking whether the canonical divisor is nef or ample. In
fact one can ask more generally what the nef cone is. The Picard group
of ${\cal A}_g^{\ast}, g\ge 2$ is generated (modulo torsion, though
for $g>2$ there is none) by two elements, namely
the $(\QQ\text{-})$ line bundle $L$ given by modular forms of weight 1
and the boundary $D$. In \cite{H2} one of us formulated the
\begin{conjecture}
A divisor $aL-bD$ on ${\cal A}_g^{\ast}$ is nef if and only if $b\ge 0$ and
$a-12 b\ge 0$.
\end{conjecture}

It is easy to see that these conditions are necessary. They are also
known to be sufficient for $g=2$ and $3$ (see below). Since the
natural quotient map ${\cal A}^{\ast}_g(n)\rightarrow {\cal
A}^{\ast}_g$ is branched of order $n$ along the boundary, this is
equivalent to
\begin{conjecture} A divisor $aL-bD$ on
${\cal A}^{\ast}_g(n)$ is nef if and only if $b\ge 0$ and
$a-12 b/n\ge 0$.
\end{conjecture}

As we shall see below one can give a quick proof of this conjecture
for $g=2$ and $3$ using known results about $\overline{\cal M}_g$ and
the Torelli map. However this approach cannot be generalized to
higher genus since the Torelli map is then no longer surjective, nor
to other than principal polarizations. For this reason an alternative
proof was given in \cite{H2} making essential use of a result of
Weissauer \cite{We} on the existence of cusp forms of small slope
which do not vanish on a given point in Siegel space.
\begin{theorem}\label{theo5}
Let $g=2$ or $3$. Then a divisor $aL-bD$ on ${\cal A}^{\ast}_g$ is nef if and
only if $b\ge 0$ and $a-12 b\ge 0$.
\end{theorem}
\begin{Proof}
First note that the two conditions are necessary. In fact let $C$ be a curve
which is contracted under the natural map $\pi:{\cal A}^{\ast}_g\rightarrow
\bar{\cal A}_g$ onto the Satake compactification. The divisor $-D$ is
$\pi$-ample (cf. also \cite{Mu4}) and $L$ is the pull-back of a line bundle on
$\bar{\cal A}_g$. Hence $(aL-bD).C\ge 0$ implies $b\ge 0$. Let $C$ be the
closure of the locus given by split abelian varieties $E\times A'$ where $E$
is an arbitrary elliptic curve and $A'$ is a fixed abelian variety of
dimension $g-1$. Then $C$ is a rational curve with $D.C=1$ and $L.C=1/12$.
This shows that $a-12 b\ge 0$ for every nef divisor $D$.

To prove that the conditions stated are sufficient we consider the
Torelli map $t:{\cal M}_g\rightarrow{\cal A}_g$ which extends to a map
$\bar{t}:\overline{\cal M}_g\rightarrow {\cal A}^{\ast}_g$. This map
is surjective for $g=2,3$. Here $\overline{\cal M}_g$ denotes the
compactification of ${\cal M}_g$ by stable curves. It follows that
for every curve $C$ in ${\cal A}^{\ast}_g$ there exists a curve $C'$
in $\overline{\cal M}_g$ which is finite over $C$. Hence a divisor on
${\cal A}^{\ast}_g, g=2,3$ is nef if and only if this is true for its
pull-back to $\overline{\cal M}_g$. We can now use Faber's paper
\cite{Fa}. Then $\bar{t}^{\ast} L=\lambda$ where $\lambda$ is the
Hodge bundle and $\bar{t}^{\ast} D=\delta_0$. Here $\delta_0$ is the
boundary $(g=2)$, resp. the closure of the locus of genus 2 curves
with one node $(g=3)$. The result now follows from \cite{Fa} since
$a\lambda-b\delta_0$ is nef on $\overline{\cal M}_g$, $g=2,3$ if $b\ge
0$ and $a-12 b\ge 0$.
\end{Proof}

\begin{cor}\label{theo6}
  The canonical divisor on ${\cal A}^{\ast}_2(n)$ is nef but not ample for
  $n=4$ and ample for $n\ge 5$. In particular ${\cal A}^{\ast}_2(4)$ is a
  minimal model and ${\cal A}^{\ast}_2(n)$ is a canonical model for $n\ge
  5$.
\end{cor}
This was first observed, though not proved in detail, by Borisov in an
early version of \cite{Bori}.

\begin{cor}\label{theo7}
  The canonical divisor on ${\cal A}^{\ast}_3(n)$ is nef but not ample for
  $n=3$ and ample for $n\ge 4$. In particular ${\cal A}^{\ast}_3(3)$ is a
  minimal model and ${\cal A}^{\ast}_3(n)$ is a canonical model for $n\ge
  4$
\end{cor}

\noindent{\em Proof of the corollaries.} Nefness or ampleness of $K$
follows immediately from Theorem \ref{theo5} since
$$
(g+1)-\frac{12}{n}\ge 0\Leftrightarrow\left\{
\begin{array}{ccc}
n\ge 4 & \mbox{ if } & g=2\\
n\ge 3 & \mbox{ if } & g=3.
\end{array}
\right.
$$
To see that $K$ is not ample on ${\cal A}^{\ast}_2(4)$ nor on ${\cal
A}^{\ast}_3(3)$ we can again use the curves $C$ coming from products
$E\times A'$ where $A'$ is a fixed abelian variety of dimension
$g-1$. For these curves $K.C=0.$\hfill$\Box$

In \cite{H3} the methods of \cite{H2} were used to prove ampleness of $K$ in
the case of $(1,p)$-polarized abel\-ian surfaces with a canonical level
structure and a
level-$n$ structure, for $p$ prime and $n \geq 5$,
provided $p$ does not divide $n$.

Finally we want to mention some results concerning the Chow ring of ${\cal
A}^{\ast}_g$. The Chow groups considered here are defined as the invariant
part of the Chow ring of ${\cal A}^{\ast}_g(n)$. The Chow ring of
$\overline{\cal M}_2$ was computed by Mumford \cite{Mu5}. This gives also
the Chow ring of ${\cal A}^{\ast}_2$, which was also calculated by a
different method by van der Geer in~\cite{vdG3}.

\begin{theorem}[\cite{Mu5},\cite{vdG3}]\label{theo8}
Let $\lambda_1=\lambda$ and $\lambda_2$ be the tautological classes on
${\cal A}^{\ast}_2$. Let $\sigma_1$ be the class of the boundary. Then
$$
\on{CH}_{\QQ}({\calA}^{\ast}_2)\cong\QQ[\lambda_1, \lambda_2, \sigma_1]/I
$$
where $I$ is the ideal generated by the relations
$$
\begin{array}{c}
(1+\lambda_1 + \lambda_2) (1-\lambda_1+\lambda_2)=1, \\
\lambda_2  \sigma_1 = 0, \\
\sigma_1^2 = 22 \sigma_1 \lambda_1-12 0\lambda_1^2.
\end{array}
$$
The ranks of the Chow groups are $1, 2, 2, 1$.
\end{theorem}
Van der Geer also computed the Chow ring of ${\cal A}^{\ast}_3$.
\begin{theorem}[\cite{vdG3}]\label{theo9}
Let $\lambda_1, \lambda_2, \lambda_3$ be the tautological classes in
${\cal A}^{\ast}_3$ and $\sigma_1, \sigma_2$ be the first and second
symmetric functions in the boundary divisors (viewed as an invariant
class on ${\cal A}^{\ast}_g(n)$). Then
$$
\on{CH}_{\QQ}({\cal A}^{\ast}_3)\cong \QQ[\lambda_1, \lambda_2,
\lambda_3,\sigma_1, \sigma_2]/ J
$$
where $J$ is the ideal generated by the relations
$$\begin{array}{lcl}
\multicolumn{3}{l}{(1+\lambda_1+\lambda_2+\lambda_3) (1-\lambda_1+
\lambda_2-\lambda_3)=1},\\
\lambda_3\sigma_1 & = & \lambda_3 \sigma_2=\lambda_1^2 \sigma_2=0,\\
\sigma_1^3 & = & 2016\lambda_3 - 4 \lambda_1^2 \sigma_1 -
24\lambda_1 \sigma_2+\frac{11}{3} \sigma_2 \sigma_1,\\
\sigma_2^2 & = & 360 \lambda_1^3  \sigma_1 -
45\lambda_1^2 \sigma_1^2+15\lambda_1 \sigma_2 \sigma_1.
\end{array}
$$
The ranks of the Chow groups are $1, 2, 4, 6, 4, 2, 1.$
\end{theorem}
\begin{Proof}
See \cite{vdG3}. The proof uses in an essential way the description of
the Voronoi compactification ${\cal A}^{\ast}_3$ given by Nakamura
\cite{Nak1} and Tsushima \cite{Ts}.
\end{Proof}

\subsection{Fundamental groups}\label{II.3}
\setcounter{theorem}{0}
The fundamental group of a smooth projective model
$\widetilde\cA(\Gamma)$ of $\cA(\Gamma)$ is independent of the
specific model chosen. We assume in this section that $g\ge 2$, so
that the dimension of $\cA(\Gamma)$ is at least~$3$.

The first results about the fundamental group of
$\widetilde\cA(\Gamma)$ were obtained by Heidrich and Kn\"oller
\cite{HK}, \cite{Kn} and concern the principal congruence subgroups
$\Gamma(n)\subset\Sp(2g,\ZZ)$. They proved the following result.

\begin{theorem}[\cite{HK},\cite{Kn}]\label{theo10}
If $n\ge 3$ or if $n=g=2$ then $\widetilde\cA_g(n)$ is simply-connected.
\end{theorem}

As an immediate corollary (first explicitly pointed out by
Heidrich-Riske) one has

\begin{cor}[\cite{HR}]\label{theo11}
If $\Gamma$ is an arithmetic subgroup of $\on{Sp}(2g,\QQ)$, then
the fundamental group of $\widetilde\cA(\Gamma)$ is finite.
\end{cor}

Corollary~\ref{theo11} follows from Theorem~\ref{theo10} because any
subgroup of $\on{Sp}(2g,\ZZ)$ of finite index contains a principal
congruence subgroup of some level.

\begin{Proof}The proof of Theorem~\ref{theo10} uses the fact that there is,
  up to the action of the group $\on{Sp}(2g,{\ZZ}_n)$,
  only one codimension~$1$ boundary component~$F$ in the Igusa
  compactification~$\cA_g^*(n)$. Suppose for simplicity that $n\ge 4$, so
  that $\Gamma(n)$ is neat. A small loop passing around this component can
  be identified with a loop in the fibre $\CC^*$ of $\cX(F)$ and hence with
  a generator $u_F$ of the $1$-dimensional lattice $P'(F)$. This loop
  determines an element $\gamma_F$, usually non-trivial, of
  $\pi_1\big(\cA_g(n)\big)$ (which is simply $\Gamma(n)$, since $\Gamma(n)$
  is torsion-free and hence acts freely on $\HH_g$). The element $\gamma_F$
  is in the kernel of the map
  $\pi_1\big(\cA_g(n)\big)\to\pi_1\big(\widetilde\cA_g(n)\big)$, so $u_F$
  is in the kernel of $\Gamma(n)\to\pi_1\big(\widetilde\cA_g(n)\big)$. But
  it turns out that the normalizer of $P'(F)$ in $\Gamma(n)$ is the whole
  of $\Gamma(n)$, as was shown by Mennicke \cite{Me} by a direct
  calculation.\end{Proof}
We (the authors of the present article) applied this method in
\cite{HS2} to the case of $\cA^{\lev}_{1,p}$ for $p\ge 5$
prime, where there are many codimension~$1$ boundary components. A
minor extra complication is the presence of some singularities in
$\Gamma\backslash\HH_2$, but they are easily dealt with. In \cite{S1}
one of us also considered the case of $\cA_{1,p}$. We found the
following simple result.

\begin{theorem}[\cite{HS2},\cite{S1}]\label{theo12}
If $p\ge 5$ is prime then $\tilde\cA^{\lev}_{1,p}$ and
$\tilde\cA_{1,p}$ are both simply-connected.
\end{theorem}

In some other cases one knows that $\widetilde\cA(\Gamma)$ is rational
and hence simply-connected. In all these cases, as F.~Campana pointed
out, it follows that the Satake compactification, and any other normal
model, is also simply-connected.

By a more systematic use of these ideas, one of us \cite{S1} gave a
more general result, valid in fact for all locally symmetric varieties
over~$\CC$. From it several results about Siegel modular varieties can
be easily deduced, of which Theorem~\ref{theo13} below is the most
striking.

\begin{theorem}[\cite{S1}]\label{theo13}
For any finite group $G$ there exists a $g\ge2$ and an arithmetic
subgroup $\Gamma\subset\Sp(2g,\ZZ)$ such that
$\pi_1\big(\widetilde\cA(\Gamma))\cong G$.
\end{theorem}

\begin{Proof}
We choose an $l\ge 4$ and a faithful representation
$\rho:G\to\Sp(2g,\FF_p)$ for some prime $p$ not dividing $2l|G|$. The
reduction mod~$p$ map $\phi_p:\Gamma(l)\to\Sp(2g,\FF_p)$ is surjective
and we take $\Gamma=\phi_p^{-1}\big(\rho(G)\big)$. As this is a
subgroup of $\Gamma(l)$ it is neat, and under these circumstances the
fundamental group of the corresponding smooth compactification of
$\cA(\Gamma)$ is $\Gamma/\Upsilon$, where $\Upsilon$ is a certain
subgroup of $\Gamma$ generated by unipotent elements (each unipotent
element corresponds to a loop around a boundary component). From this
it follows that $\Upsilon\subset\Ker\phi_p=\Gamma(pl)$. Then from
Theorem~\ref{theo10} applied to level $pl$ it follows that
$\Upsilon=\Gamma(pl)$ and hence that the fundamental group is
$\Gamma/\Gamma(pl)\cong G$.
\end{Proof}

For $G=D_8$ we may take $g=2$; in particular, the fundamental group
of a smooth projective model of a Siegel modular threefold need not
be abelian. Apart from the slightly artificial examples which
constitute Theorem~\ref{theo13}, it is also shown in \cite{S1} that a
smooth model of the double cover $\tilde \cN_5$ of Nieto's threefold
$\cN_5$ has fundamental group $\ZZ_2\times\ZZ_2$. The space
$\tilde\cN_5$ will be discussed in Section~IV below: it is birational
with the moduli space of abelian surfaces with a polarization of
type~$(1,3)$ and a level-$2$ structure.

\renewcommand{\thetheorem}{\Roman{section}.\arabic{subsection}.\arabic{theorem}}
\addtocounter{theorem}{1}
\section{Abelian surfaces}\label{III}
In the case of abelian surfaces the moduli spaces ${\cal A}_{1,t}$ and
${\cal A}_{1,t}^{\lev}$ of abelian surfaces with a
$(1,t)$-polarization, resp. with a $(1,t)$-polarization and a
canonical level structure were investigated by a number of
authors. One of the starting points for this development was the paper
by Horrocks and Mumford \cite{HM} which established a connection
between the Horrocks-Mumford bundle on $\PP^4$ and the moduli space
${\cal A}_{1,5}^{\lev}$.

\subsection{The lifting method}\label{III.1}
\setcounter{theorem}{0}
Using a version of Maa\ss{} lifting Gritsenko has proved the existence
of a weight 3 cusp forms for almost all values of $t$. Before we can
describe his lifting result recall the {\it paramodular group}
$\on{Sp}(\Lambda,\ZZ)$ where
$$
\Lambda=\left(
\begin{array}{cc}
0  &  E\\
-E &  0
\end{array}
\right), \, E=\left(
\begin{array}{cc}
1 & 0\\
0 & t
\end{array}
\right)
$$
for some integer $t\ge 1$, with respect to a basis
$(e_1,e_2,e_3,e_4)$. This group is conjugate to the (rational)
paramodular group
$$
\Gamma_{1,t}=R^{-1} \on{Sp}(\Lambda,\ZZ) R,\qquad R={\left(
\begin{smallmatrix}
1&&&\\
&1&&\\
&&1&\\
&&&t
\end{smallmatrix}
\right)}.
$$
It is straightforward to check that
$$
\Gamma_{1,t} = \left\{
g\in \on{Sp}(4,\QQ); \, g\in \left(
\begin{array}{cccc}
\ZZ  & \ZZ & \ZZ  & t\ZZ\\
t\ZZ & \ZZ & t\ZZ & t\ZZ\\
\ZZ  & \ZZ & \ZZ  & t\ZZ\\
\ZZ  & t^{-1}\ZZ  & \ZZ & \ZZ
\end{array}
\right)\right\}.
$$
Then ${\cal A}_{1,t}=\Gamma_{1,t}\backslash \HH_2$ is the moduli space of
$(1,t)$-polarized abelian surfaces. In this chapter we shall denote the
elements of $\HH_2$ by
$$
\tau=\left(
\begin{array}{cc}
\tau_1 & \tau_2\\
\tau_2 & \tau_3
\end{array}
\right)\in \HH_2.
$$
The Tits building of $\Gamma_{1,t}$, and hence the combinatorial
structure of the boundary components of the Satake or the Voronoi
(Igusa) compactification of ${\cal A}_{1,t}$ are known, at least if
$t$ is square free: see \cite{FrS}, where Tits buildings for some other
groups are also calculated. There are exactly $\mu (t)$ corank $1$ boundary
components, where $\mu (t)$ denotes the number of prime divisors of
$t$ \cite[Folgerung 2.4]{Gr1}. If $t$ is square free, then there
exists exactly one corank $2$ boundary component \cite[Satz
4.7]{Fr}. In particular, if $t>1$ is a prime number then there exist
two corank $1$ boundary components and one corank $2$ boundary
component. These boundary components belong to the isotropic subspaces
spanned by $e_3$ and $e_4$, resp. by $e_3\wedge e_4$. In terms of the
Siegel space the two corank $1$ boundary components correspond to
$\tau_1\rightarrow i\infty$ and $\tau_3\rightarrow i
\infty$. For $t=1$ these two components are equivalent under the group
$\Gamma_{1,1}=\on{Sp}(4,\ZZ)$.

Gritsenko's construction of cusp forms uses a version of Maa\ss{}
lifting. In order to explain this, we first have to recall the
definition of {\em Jacobi forms}. Here we restrict ourselves to the
case of $\Gamma_{1,1}=\on{Sp}(4,\ZZ)$. The stabilizer of $\QQ e_4$ in
$\on{Sp}(4,\ZZ)$ has the structure
$$
P(e_4)\cong \on{SL}(2,\ZZ)\ltimes H(\ZZ)
$$
where $\SL(2,\ZZ)$ is identified with
$$
\left\{
\left(
\begin{array}{cccc}
a & 0 & b & 0\\
0 & 1 & 0 & 0\\
c & 0 & d & 0\\
0 & 0 & 0 & 1
\end{array}\right);
\left(
\begin{array}
{cc}
a & b\\
c & d
\end{array}\right)
\in \on{SL}(2,\ZZ)
\right\}
$$
and
$$
H(\ZZ)=\left\{
\left(
\begin{array}{cccc}
1 & 0 & 0 & \mu\\
\lambda & 1 & \mu & r\\
0 & 0 & 1 & -\lambda\\
0 & 0 & 0 & 1
\end{array}\right);
\lambda, \mu, r \in \ZZ
\right\}
$$
is the {\it integral Heisenberg group}.

Every modular form $F\in M_k(\on{Sp}(4, \ZZ))$ of weight $k$ with respect
to $\on{Sp}(4, \ZZ)$ has a Fourier extension with respect to $\tau_3$
which
is of the following form
$$
F(\tau)=\sum\limits_{m\ge 0}f_m(\tau_1, \tau_2) e^{2\pi im\tau_3}.
$$
The same is true for modular forms with respect to $\Gamma_{1,t}$, the
only difference is that the factor $\exp(2\pi im \tau_3)$ has to be
replaced by $\exp(2\pi im t \tau_3)$. The coefficients $f_m(\tau_1, \tau_2)$
are examples of {\em Jacobi forms}. Formally Jacobi forms are defined as
follows:
\begin{definition}
  A {\em Jacobi form} of {\em index} $m$ and {\em weight} $k$ is a
  holomorphic function
$$
\Phi=\Phi(\tau , z):\HH_1\times \CC\rightarrow \CC
$$
which has the following properties:
\begin{enumerate}
\item[(1)]
It has the {\em transformation behaviour}
\begin{enumerate}
\item[(a)]
$\Phi\left(\frac{a\tau + b}{c\tau+d}, \frac {z}{c\tau+d} \right)=
(c\tau+d)^k e^{\frac{2\pi i c m z^2}{c\tau+d}}\Phi(\tau, z),\left(
\begin{array}{cc} a & b\\ c & d\end{array}\right) \in \on{SL}(2,\ZZ)$
\item[(b)]
$\Phi(\tau, z+\lambda\tau+\mu)=e^{-2\pi i m (\lambda^2\tau+2\lambda z)}
\Phi(\tau,z),\quad \lambda, \mu \in \ZZ$.
\end{enumerate}
\item[(2)]
It has a Fourier expansion
$$
\Phi(\tau, z)=\sum\limits_{\begin{array}{l}
{{}^{n,l\in\ZZ, n\ge 0}}\\[-2mm]
{{}^{4 n m\ge l^2}}\end{array}} f(n,l)e^{2\pi i(n\tau+lz)}.
$$
\end{enumerate}
A Jacobi form is called a  {\em cusp form} if one has strict inequality $4n
m>l^2 $ in the Fourier expansion.
\end{definition}

Note that for $z=0$ the transformation behaviour described by (1)(a)
is exactly that of a modular form. For fixed $\tau$ the transformation
law (1)(b) is, up to a factor $2$ in the exponent, the transformation
law for theta functions. One can also summarize (1)(a) and (1)(b) by
saying that $\Phi=\Phi(\tau, z)$ is a modular form with respect to the
{\em Jacobi group} $\on{SL}(2,\ZZ)\ltimes H(\ZZ)$. (Very roughly,
Jacobi forms can be thought of as sections of a suitable $\QQ$-line
bundle over the universal elliptic curve, which doesn't actually
exist.) The Jacobi forms of weight $k$ and index $m$ form a vector
space $J_{k,m}$ of finite dimension. The standard reference for Jacobi
forms is the book by Eichler and Zagier \cite{EZ}.

As we have said before, Jacobi forms arise naturally as coefficients
in the Fourier expansion of modular forms. These coefficients are
functions, or more precisely sections of a suitable line bundle, on a
boundary component of the Siegel modular threefold. The idea of
lifting is to reverse this process. Starting with a Jacobi form one
wants to construct a Siegel modular form where this Jacobi form
appears as a Fourier coefficient. This idea goes back to
Maa\ss{}~\cite{Ma2} and has in recent years been refined in several
ways by Gritsenko, Borcherds and others: see e.g.\cite{Gr1},
\cite{Gr3}, \cite{GrN} and \cite{Borc}. The following lifting result
is due to Gritsenko.
\begin{theorem}[\cite{Gr1}]\label{theo14}
There is a lifting, i.e. an embedding
$$
\operatorname{Lift }: J_{k,t}\longrightarrow M_k (\Gamma_{1,t})
$$
of the space of Jacobi forms of weight $k$ and index $t$ into the space of
modular forms of weight $k$ with respect to the paramodular group
$\Gamma_{1,t}$. The lifting of a Jacobi cusp form is again a cusp form.
\end{theorem}

\begin{Proof}
For details see \cite[Hauptsatz 2.1]{Gr1} or \cite[Theorem 3]{Gr2}. For a
Jacobi form $\Phi=\Phi(\tau z)
$ with Fourier expansion
$$
\Phi(\tau, z)=\sum\limits_{\begin{array}{l}
{{}^{n,l\in\ZZ}}\\[-2mm]
{{}^{4 n t\ge l^2}}\end{array}} f(n,l)e^{2\pi i(n\tau+lz)}
$$
the lift can be written down explicitly as
$$
\on{Lift} \Phi(\tau)=\sum\limits_{4 t m n \ge l^2} \sum\limits_{a|(n,
l, m)} a^{k-1} f\left(\frac{n m}{a^2}, \frac{l}{a}\right) e^{2\pi
i(n\tau_1+l\tau_2+m t\tau_3)}.
$$
\end{Proof}

\noindent Since one knows dimension formulae for Jacobi cusp forms one
obtains in this way lower bounds for the dimension of the space of
modular forms and cusp forms with respect to the paramodular
group. Using this together with Freitag's extension theorem it is then
easy to obtain the following corollaries.
\begin{cor}\label{theo15}
  Let $p_g(t)$ be the geometric genus of a smooth projective model of the
  moduli space ${\cal A}_{1,t}$ of $(1,t)$-polarized abelian surfaces. Then
$$
p_g(t)\ge\sum\limits^{t-1}_{j=1}\left(
\{2j+2\}_{12}-\left\lfloor\frac{j^2}{12}\right\rfloor\right)
$$
where
$$
\{m\}_{12}=\left\{
\begin{array}{lll}
\left\lfloor\frac m{12}\right\rfloor     & \mbox{if} & m\not\equiv 2
\mod 12 \\[1mm]
\left\lfloor\frac m{12}\right\rfloor  -1 & \mbox{if} & m\equiv 2 \mod 12
\end{array}
\right .
$$
and $\lfloor x\rfloor$ denotes the integer part of $x$.
\end{cor}

This corollary also implies that $p_g(t)$ goes to infinity as $t$ goes to
infinity.
\begin{cor}\label{theo16}
  The Kodaira dimension of ${\cal A}_{1,t}$ is non-negative if $t\ge 13$
  and $t\neq 14, 15, 16, 18, 20, 24, 30$ and $36$. In particular these
  spaces are not unirational.
\end{cor}

\begin{cor}\label{theo17}
  The Kodaira dimension of ${\cal A}_{1,t}$ is positive if $t\ge 29$ and
  $t\neq 30, 32, 35, 36, 40, 42, 48$ and $60$.
\end{cor}

On the other hand one knows that ${\cal A}_{1,t}$ is rational or
unirational for small values of $t$. We have already mentioned that Igusa
proved rationality of ${\cal A}_{1,1}={\cal A}_2$ in \cite{I1}.
Rationality of ${\cal A}_{1,2}$ and ${\cal A}_{1,3}$ was proved by
Birkenhake and Lange \cite{BL}. Birkenhake, Lange and van Straten
\cite{BLvS} also showed that ${\cal A}_{1,4}$ is unirational. It is a
consequence of the work of Horrocks and Mumford \cite{HM} that ${\cal
A}_{1,5}^{\lev}$ is rational. The variety ${\cal A}_{1,7}^{\lev}$ is
birational to a Fano variety of type $V_{22}$ \cite{MS} and hence also
rational. The following result of Gross and Popescu was stated in \cite{GP1}
and is proved in the series of papers \cite{GP1}--\cite{GP3}.

\begin{theorem}[\cite{GP1},\cite{GP1a},\cite{GP2},\cite{GP3}]\label{theo18}
${\cal A}_{1,t}^{\lev}$ is rational for $6\le
t\le 10$ and $t=12$ and unirational, but not rational, for $t=11$.
Moreover the variety ${\cal A}_{1,t}$ is unirational for $t=14, 16,
18$ and $20$.
\end{theorem}

We shall return to some of the projective models of the modular
varieties ${\cal A}_{1,t}$ in chapter~V. Altogether this gives a
fairly complete picture as regards the question which of the spaces
${\cal A}_{1,t}$ can be rational or unirational. In fact there are
only very few open cases.

\begin{problem}
Determine whether the spaces ${\cal A}_{1,t}$ for $t=15, 24, 30$ or
$36$ are unirational.
\end{problem}

\subsection{General type results for moduli spaces of abelian
surfaces}\label{III.2}
\setcounter{theorem}{0}
In the case of moduli spaces of abelian surfaces there are a number of
concrete bounds which guarantee that the moduli spaces ${\cal
A}_{1,t}$, resp. ${\cal A}_{1,t}^{\lev}$ are of
general type. Here we collect the known results and comment on the
different approaches which enable one to prove these theorems.

\begin{theorem}[\cite{HS1},\cite{GrH1}]\label{theo19}
Let $p$ be a prime number. The moduli spaces ${\cal
A}_{1,p}^{\lev}$ are of general type if $p\ge 37$.
\end{theorem}

\begin{Proof}This theorem was first proved in \cite{HS1} for $p\ge 41$ and
  was improved in \cite{GrH1} to $p=37$. The two methods of proof differ in
  one important point.  In \cite{HS1} we first estimate how the dimension
  of the space of cusp forms grows with the weight $k$ and find that
\begin{equation}
\dim S_{3k} \left(\Gamma_{1,p}^{\lev}\right) =
\frac{p(p^4-1)}{640} k^3+O(k^2).
\end{equation}
These cusp forms give rise to $k$-fold differential forms on $\opn{\cal
A}_{1,p}^{\lev}$ and we have two types of obstruction to extending them to
a smooth projective model of ${\cal A}_{1,p}^{\lev}$: one comes from the
boundary and the other arises from the elliptic fixed points. To calculate
the number of obstructions from the boundary we used the description of the
boundary of the Igusa compactification (which is equal to the Voronoi
decomposition) given in \cite{HKW2}. We found that the number of
obstructions to extending $k$-fold differentials is bounded by
\begin{equation}
H_B(p,k)=\frac{(p^2-1)}{144}(9p^2+2p+11) k^3+O(k^2).
\end{equation}
The singularities of the moduli spaces ${\cal A}_{1,p}^{\lev}$ and of the
Igusa compactification were computed in \cite{HKW1}. This allowed us to
calculate the obstructions arising from the fixed points of the action of
the group $\Gamma_{1,p}^{\lev}$.  The result is that the number of these
obstructions is bounded by
\begin{equation}
H_S(p,k)=\frac{1}{12}(p^2-1)(\frac{7}{18} p-1) k^3+O(k^2).
\end{equation}
The result then follows from comparing the leading terms of
(6) and (7) with that of (5).

The approach in \cite{GrH1} is different. The crucial point is to use
Gritsenko's lifting result to produce non-zero cusp forms of weight $2$.
The first prime where this works is $p=37$, but it also works for all
primes $p>71$. Let $G$ be a non-trivial modular form of weight $2$ with
respect to $\Gamma_{1,37}$. Then we can consider the subspace
$$
V_k = G^k M_k \left(\Gamma_{1,37}^{\lev}\right)
\subset M_{3
k} \left(\Gamma_{1,37}^{\lev}\right).
$$
The crucial point is that the elements of $V_k$ vanish by construction
to order $k$ on the boundary. This ensures that the extension to the
boundary imposes no further conditions. The only possible obstructions are
those coming from the elliptic fixed points. These obstructions were
computed above. A comparison of the leading terms again gives the
result.\end{Proof}

The second method described above was also used in the proof of the following
two results.

\begin{theorem}[\cite{OG},\cite{GrS}]\label{theo20}
The moduli space ${\cal A}_{1,p^2}^{\lev}$ is of general type for
every prime $p\ge 11$.
\end{theorem}

This was proved in \cite{GrS} and improves a result of O'Grady
\cite{OG} who had shown this for $p\ge 17$. The crucial point in
\cite{GrS} is that, because of the square $p^2$, there is a covering
${\cal A}_{1,p^2}\rightarrow {\cal A}_{1,1}$. The proof in \cite{GrS}
then also uses the existence of a weight $2$ cusp form with respect to
the group $\Gamma_{1,p^2}$ for $p\ge 11$. The only obstructions which
have to be computed explicitly are those coming from the elliptic
fixed points. The essential ingredient in O'Grady's proof is the
existence of a map from a partial desingularization of a toroidal
compactification to the space $\overline{\cal M}_2$ of semi-stable
genus $2$ curves.

A further result in this direction is

\begin{theorem}[\cite{S2}]\label{theo21}
The moduli spaces ${\cal A}_{1,p}$ are of general type for all primes
$p\ge 173$.
\end{theorem}

It is important to remark that in this case there is no natural map from
${\cal A}_{1,p}$ to the moduli space ${\cal A}_{1,1}={\cal A}_2$ of
principally polarized abelian surfaces. A crucial ingredient in the proof
of the above theorem is the calculation of the singularities of the spaces
${\cal A}_{1,p}$ which was achieved by Brasch \cite{Br}. Another recent
result is

\begin{theorem}[\cite{H3}]\label{theo21a}
The moduli spaces of $(1,d)$-polarized abel\-ian
surfaces with a full level-$n$ structure are of general type for all pairs
$(d,n)$ with $(d,n)=1$ and $n \geq 4$.
\end{theorem}

A general result due to L.~Borisov is

\begin{theorem}[\cite{Bori}]\label{theo22}
There are only finitely many subgroups $H$ of $\on{Sp}(4,\ZZ)$ such that
${\cal A}(H)$ is not of general type.
\end{theorem}

Note that this result applies to the groups $\Gamma_{1,p}^{\lev}$ and
$\Gamma_{1,p^2}$ which are both conjugate to subgroups of $\on{Sp}(4,
\ZZ)$, but does not apply to the groups $\Gamma_{1,p}$, which are not. (At
least for $p\geq 7$: the subgroup of $\CC^*$ generated by the eigenvalues
of non-torsion elements of $\Gamma_{1,p}$ contains $p$th roots of unity, as
was shown by Brasch in~\cite{Br}, but the corresponding group for
$\Sp(4,\ZZ)$ has only $2$- and $3$-torsion.)

We shall give a rough outline of the proof of this result. For details the
reader is referred to \cite{Bori}. We shall mostly comment on the geometric
aspects of the proof. Every subgroup $H$ in $\on{Sp}(4,\ZZ)$ contains a
principal congruence subgroup $\Gamma(n)$. The first reduction is the
observation that it is sufficient to consider only subgroups $H$ which contain
a principal congruence subgroup $\Gamma(p^t)$ for some prime $p$. This
is essentially a group theoretic argument using the fact that the finite
group $\on{Sp}(4,{\ZZ}_p)$ is simple for all primes $p\geq3$. Let us now
assume that $H$ contains $\Gamma(n)$ (we assume $n\geq5$).
This implies that there is a finite morphism
${\cal A}_2(n) \to {\cal A}(H)$. The idea is to show that
for almost all groups $H$ there are sufficiently many pluricanonical
forms on the Igusa (Voronoi) compactification $X={\cal A}_2^{\ast}(n)$
which descend to a smooth projective model of ${\cal A}(H)$. For this it is
crucial to get a hold on the possible singularities of the quotient $Y$.
We have already observed in Corollary \ref{theo7} that the
canonical divisor on $X$ is ample for $n\geq 5$. The finite group
$\bar{H}=\Gamma_2(n)/H$ acts on $X$ and the
quotient $Y=\bar{H}\backslash X$ is a
(in general singular) projective model of ${\cal A}(H)$. Since $X$ is smooth
and $H$ is finite, the variety $Y$ is normal and has log-terminal
singularities, i.e. if $\pi:Z\rightarrow Y$ is a desingularization whose
exceptional divisor $E=\sum\limits_i E_i$ has simple normal crossing, then
$$
K_Z=\pi^{\ast} K_Y+\sum\limits_i(-1+\delta_i) E_i\quad\mbox{ with }
\delta_i>0.
$$
Choose $\delta>0$ such that $-1+\delta$ is the minimal discrepancy. By $L_X$,
resp.  $L_Y$ we denote the $\QQ$-line bundle whose sections are modular
forms of weight $1$. Then $L_X=\mu^{\ast}L_Y$ where $\mu:X\rightarrow Y$ is
the quotient map.

The next reduction is that it suffices to construct a non-trivial
section $s\in H^0(m(K_Y-L_Y))$ such that $s_y \in {\cal
O}_Y \left(m(K_Y-L_Y)\gothm_y^{m(1-\delta)}\right)$ for all $y\in Y$
where $Y$ has a non-canonical singularity. This is enough because
$\pi^{\ast}(sH^0(m L_Y))\subset H^0(m K_Z)$ and the
dimension of the space $H^0(m L_Y))$ grows as $m^3$.

The idea is to construct $s$ as a suitable $\bar{H}$-invariant section
$$
s\in H^0\big(\mu^*(m(K_Y-L_Y))\big)^{\bar H}
$$
satisfying vanishing conditions
at the branch locus of the finite map $\mu:X\rightarrow Y$. For this
one has to understand the geometry of the quotient map~$\mu$. First of
all one has branching along the boundary $D=\sum D_i$ of
$X$. We also have to look at the Humbert surfaces
$$
{\cal H}_1=\left\{ \tau=\left(
\begin{array}{cc}
\tau_1  &  0\\
0       &  \tau_3
\end{array}\right); \tau_1, \tau_3\in \HH_1\right\}=\mbox{ Fix }
\left(
\begin{array}{cccc}
1 &&&\\
&-1&&\\
&&1 &\\
&&&-1
\end{array}
\right)
$$
and
$$
{\cal H}_4=\left\{ \tau=\left(
\begin{array}{cc}
\tau_1  &  \tau_2\\
\tau_2  &  \tau_3
\end{array}\right); \tau_1 = \tau_3\right\}=\mbox{ Fix }
\left(
\begin{array}{cccc}
0 &1&&\\
1 &0&&\\
&&0&1\\
&&1&0
\end{array}
\right).
$$
Let
$$
{\cal F}=\bigcup\limits_{g\in \on{Sp}(4,\ZZ)} g({\cal H}_1)\quad,\quad
{\cal
G}=\bigcup\limits_{g\in \on{Sp}(4,\ZZ)} g({\cal H}_2)
$$
and let
$$
F=\overline{\pi({\cal F})},\quad  G=\overline{\pi({\cal G})}
$$
where $\pi:\HH_2\rightarrow \Gamma(n)\backslash\HH_2\subset X$ is the
quotient map. One can then show that the branching divisor of the map
${\cal A}(\Gamma_2(n))\rightarrow {\cal A}(H)$ is contained in $F\cup G$
and that all singularities in ${\cal A}(H)$ which lie outside $\mu (F\cup
G)$ are canonical. Moreover the stabilizer subgroups in $\on{Sp}(4, \ZZ)$
of points in ${\cal F} \cup {\cal G}$ are solvable groups of bounded order.
Let $F=\sum F_i$ and $G=\sum G_i$ be the decomposition of the surfaces $F$
and $G$ into irreducible components. We denote by $d_i, f_i$ and $g_i$ the
ramification order of the quotient map $\mu:X\rightarrow Y$ along $D_i,
F_i$ and $G_i$. The numbers $f_i$ and $g_i$ are equal to $1$ or $2$.
One has
\begin{eqnarray*}
\mu^{\ast}(m(K_Y-L_Y))&=&m(K_X-L_X) - \sum\limits_i m(d_i-1)D_i-\sum\limits_i
m(f_i-1) F_i\\
&&\mbox{  }- \sum\limits_i m(g_i-1) G_i.
\end{eqnarray*}

Recall that the finite group $\bar{H}$ is a subgroup of the group
$\bar{G}=\Gamma/\Gamma(n)=\on{Sp}(4,\ZZ_n)$.
The crucial point in Borisov's argument
is to show, roughly speaking, that the index $[\bar{G}:\bar{H}]$
can be bounded from above in terms of the singularities of $Y$. There
are several such types of bounds depending on whether one considers
points on the branch locus or on one or more boundary components.
We first use this
bound for the points on $X$ which lie on $3$ boundary divisors.
Using this and the fact that $Y$ has only finite quotient
singularities one obtains
the following further reduction: if $R$ is the ramification divisor of the
map $\mu: X \to Y$, then it is enough to construct a non-zero
section in $H^0(m(K_X-L_X-R))$ for some $m>0$ which lies in
$\gothm_x^{m k(\on{Stab}^{H}x)}$ for all points $x$ in $X$ which lie over
non-canonical points of $Y$ and which are not on the intersection
of $3$ boundary divisors. Here $k(\on{Stab}^{H}x)$ is
defined as follows. First note that $\on{Stab}^H x$ is solvable and
consider a series
$$
\{0\}=H_0\triangleleft H_1\triangleleft\ldots\triangleleft
H_t=\on{Stab}^H x
$$
with $H_i/H_{i-1}$ abelian of exponent $k_i$. Take
$k'=k_1\cdot\ldots\cdot k_t$. Then $k(\on{Stab}^{H}x)$ is the minimum over
all $k'$ which are obtained in this way.
To obtain an invariant section one can then
take the product with respect to the action of the finite group $\bar H$.
Now recall that all non-canonical points on ${\cal A}(H)$ lie in
$\mu (F \cup G)$. The subgroup $Z\on{Stab}^H D_i$ of
$\on{Stab}^H D_i$ which acts trivially on $D_i$ is cyclic of order $d_i$.
Moreover if $x$ lies on exactly one boundary divisor
of $X$ then the order of the group $\on{Stab}^H x/Z \on{Stab}^H D_i$ is
bounded by $6$ and if $x$ lies on exactly $2$ boundary divisors, then
the order of this group is bounded by $4$. Using this one can show that
there is a constant $c$ (independent of $H$) such that
it is sufficient to construct a non-zero section in $m(K_X-L_X-cR)$ for some
positive $m$.
By results of Yamazaki \cite{Ya} the
divisor $mK_X - 2 m L_X$ is effective. It is,
therefore, sufficient to prove the existence of a non-zero section in
$m(K_X-2cR)$. The latter equals
$$
mK_X- 2c\sum\limits_i m(d_i-1)D_i-2c\sum\limits_i
m(f_i-1) F_i)- 2c\sum\limits_i m(g_i-1) G_i.
$$
We shall now restrict ourselves to
obstructions coming from components $F_i$; the obstructions coming from
$G_i$, $D_i$ can be treated
similarly. Since $h^0(mK_X)> c_1 n^{10}m^3$ for some $c_1 > 0$, $m \gg 0$ one
has to prove the following result: Let $\eps >0$. Then for all but finitely
many subgroups $H$ one has
$$
  \sum_{f_i=2} (h^0 (mK_X)- h^0(mK_X-2cm f_i F_i)) \leq \eps n^{10} m^3\quad
  \mbox{for}\quad m\gg 0
$$
and all $n$. This can finally be derived from the following boundedness result.
Let $\eps >0$ and assume that
$$
  \frac{\#\{F_i\>;\>\> f_i=2\}}{\# \{F_i\}}\geq \eps ,
$$
then the index $[\bar G : \bar H]$ is bounded by an (explicitly known) constant
depending only on $\eps$.
The proof of
this statement is group theoretic and the idea is as follows. Assume the
above inequality holds: then $H$ contains many involutions and these
generate a subgroup of $\on{Sp}(4,\ZZ)$ whose index is bounded in terms
of~$\eps$.

\subsection{Left and right neighbours}\label{III.3}
\setcounter{theorem}{0}
The paramodular group $\Gamma_{1,t}\subset \on{Sp}(4,\QQ)$ is (for
$t>1$) not a maximal discrete subgroup of the group of analytic
automorphisms of $\HH_2$. For every divisor $d \| t$ (i.e. $d|t$
and $(d,t/d)=1)$ one can choose integers $x$ and $y$ such that
$$
  xd - yt_d = 1,\quad \mbox{where }t_d =t/d.
$$
The matrix
$$
  V_d = \frac{1}{\sqrt{d}}
    \begin{pmatrix}{}
      dx & -1 & 0 & 0\\
      -yt & d  & 0 & 0\\
      0   & 0  & d & yt\\
      0   & 0  & 1 & dx
    \end{pmatrix}
$$
is an element of $\on{Sp}(4,\RR)$ and one easily checks that
$$
  V_d^2 \in \Gamma_{1,t},\qquad V_d \Gamma_{1,t} V_d^{-1} = \Gamma_{1,t}.
$$
The group generated by $\Gamma_{1,t}$ and the elements $V_d$, i.e.
$$
  \Gamma_{1,t}^\ast = \langle \Gamma_{1,t},\; V_d; \, d\| t\rangle
$$
does not depend on the choice of the integers $x,y$. It is a normal extension
of $\Gamma_{1,t}$ with
$$
  \Gamma_{1,t}^{\ast}/\Gamma_{1,t} \cong (\ZZ_2)^{\nu (t)}
$$
where $\nu (t)$ is the number of prime divisors of $t$. If $t$ is
square-free, it is known that $\Gamma_{1,t}^\ast$ is a maximal discrete
subgroup of $\on{Sp}(4,\RR)$ (see \cite{All},\cite{Gu}). The coset
$\Gamma_{1,t}V_t$ equals $\Gamma_{1,t} V_t'$ where
$$
  V_d =
    \begin{pmatrix}{}
      0 & \sqrt{t}^{-1} & 0 & 0\\
      \sqrt{t} & 0  & 0 & 0\\
      0   & 0  & 0 & \sqrt{t}\\
      0   & 0  & \sqrt{t}^{-1} & 0
    \end{pmatrix}.
$$
This generalizes the Fricke involution known from the theory of
elliptic curves. The geometric meaning of the involution $\bar V_t:
\calA_{1,t}\to \calA_{1,t}$ induced by $V_t$ is that it maps a polarized
abelian surface $(A,H)$ to its dual. A similar geometric
interpretation can also be given for the involutions $V_d$ (see
\cite[Proposition 1.6]{GrH2} and also \cite[Satz (1.11)]{Br} for the case
$d=t$). We also consider the degree $2$ extension
$$
  \Gamma_{1,t}^+ =\langle \Gamma_{1,t}, V_t \rangle
$$
of $\Gamma_{1,t}$. If $t=p^n$ for a prime number $p$, then
$\Gamma_{1,t}^+=\Gamma_{1,t}^\ast$. The groups $\Gamma_{1,t}^\ast$ and
$\Gamma_{1,t}^+$ define Siegel modular threefolds
$$
 \calA_{1,t}^\ast = \Gamma_{1,t}^\ast \setminus \HH_2, \qquad
 \calA_{1,t}^+=\Gamma_{1,t}^+ \setminus \HH_2.
$$
Since $\Gamma_{1,t}^\ast$ is a maximal discrete subgroup for $t$ square free
the space $\calA_{1,t}^\ast$ was called a \emph{minimal Siegel modular
threefold}. This should not be confused with minimal models in the sense of
Mori theory.

The paper \cite{GrH2} contains an interpretation of the varieties
$\calA_{1,t}^\ast$ and $\calA_{1,t}^+$ as moduli spaces.
We start with the spaces
$\calA_{1,t}^\ast$.

\setcounter{theorem}{0}
\begin{theorem}[\cite{GrH2}]
\begin{enumerate}
\item[\rm (i)] Let $A,A'$ be two $(1,t)$-polarized abelian surfaces
which define the same point in $\calA_{1,t}^\ast$. Then their (smooth)
Kummer surfaces $X,X'$ are isomorphic.
\item[\rm (ii)] Assume that the N\'eron-Severi group of $A$ and $A'$ is
generated by the polarization. Then the converse is also true: if $A$
and $A'$ have isomorphic Kummer surfaces, then $A$ and $A'$ define the
same point in $\calA_{1,t}^\ast$.
\end{enumerate}
\end{theorem}

The proof of this theorem is given in \cite[Theorem 1.5]{GrH2}. The crucial
ingredient is the Torelli theorem for $K3$-surfaces. The above theorem
says in particular that an abelian surface and its dual have
isomorphic Kummer surfaces. This implies a negative answer to a
problem posed by Shioda, who asked whether it was true that two
abelian surfaces whose Kummer surfaces are isomorphic are necessarily
isomorphic themselves. In view of the above result, a general
$(1,t)$-polarized surface with $t>1$ gives a counterexample: the
surface $A$ and its dual $\hat A$ have isomorphic Kummer surfaces, but
$A$ and $\hat A$ are not isomorphic as polarized abelian surfaces. If
the polarization generates the N\'eron-Severi group this implies that
$A$ and $\hat A$ are not isomorphic as algebraic surfaces. In view of
the above theorem one can interpret $\calA_{1,t}^\ast$ as the space of
Kummer surfaces associated to $(1,t)$-polarized abelian surfaces.

The space $\calA_{1,t}^+$ can be interpreted as a space of
lattice-polarized $K3$-surfaces in the sense of \cite{N3},\cite{Dol}. As
usual let $E_8$ be the even, unimodular, positive definite lattice of
rank 8. By $E_8(-1)$ we denote the lattice which arises from $E_8$ by
multiplying the form with $-1$. Let $\langle n \rangle$ be the rank 1
lattice $\ZZ l$ with the form given by $l^2=n$.

\begin{theorem}[\cite{GrH2}]
The moduli space $\calA_{1,t}^+$ is isomorphic to the moduli space of
lattice polarized $K3$-surfaces with a polarization of type $\langle
2t\rangle \oplus 2 E_8 (-1)$.
\end{theorem}

For a proof see \cite[Proposition 1.4]{GrH2}. If
$$
  L= \ZZ e_1 \oplus \ZZ e_2 \oplus \ZZ e_3 \oplus \ZZ e_4 ,
$$
then $\bigwedge\nolimits^2 L$ carries a symmetric bilinear form $(\,\, ,
\, )$ given by
$$
  x \wedge y = (x,y) e_1 \wedge e_2 \wedge e_3 \wedge e_4 \in
  \bigwedge\nolimits^4 L.
$$
If $w_t = e_1 \wedge e_3+ t e_2 \wedge e_4$, then the group
$$
  \tilde \Gamma_{1,t} =\lbrace g:\, L\to L; \,\bigwedge\nolimits^2
g(w_t)=w_t \rbrace
$$
is isomorphic to the paramodular group $\Gamma_{1,t}$. The lattice
$L_t = w_t^\perp$ has rank~$5$ and the form $(\,\, , \, )$ induces a
quadratic form of signature (3,2) on $L_t$. If $O(L_t)$ is the
orthogonal group of isometries of the lattice $L_t$, then there is a
natural homomorphism
$$
  \bigwedge\nolimits^2: \quad \Gamma_{1,t} \cong \tilde \Gamma_{1,t}
\longrightarrow
O(L_t).
$$
This homomorphism can be extended to $\Gamma_{1,t}^\ast$ and
$$
  \Gamma_{1,t}^\ast/\Gamma_{1,t}
\cong O(L^\vee_t / L_t) \cong (\ZZ_2)^{\nu(t)}
$$
where $L^\vee_t$ is the dual lattice of $L_t$. This, together with
Nikulin's theory (\cite{N2}, \cite{N3}) is the crucial ingredient in the
proof of the above theorems.

The varieties $\calA_{1,t}^+$ and $\calA_{1,t}^\ast$ are quotients of the
moduli space $\calA_{1,t}$ of $(1,t)$-polarized abelian surfaces. In
\cite{GrH3} there is an investigation into an interesting class of Galois
coverings of the spaces $\calA_{1,t}$. These coverings are called
\emph{left neighbours}, and the quotients are called \emph{right
neighbours}. To explain the coverings of $\calA_{1,t}$ which were
considered in \cite{GrH3}, we have to recall a well known result about the
commutator subgroup $\on{Sp}(2g,\ZZ)'$ of the symplectic group
$\on{Sp}(2g,\ZZ)$. Reiner \cite{Re} and Maa\ss{}~\cite{Ma1} proved that
$$
  \on{Sp}(2g,\ZZ)/ \on{Sp}(2g,\ZZ)'= \begin{cases}
    \ZZ_{12} & \mbox{ for } g=1\\
    \ZZ_2 & \mbox{ for } g=2\\
    1 & \mbox{ for }g \geq 3
  \end{cases}.
$$
The existence of a character of order 12 of $\on{Sp}(2,\ZZ)=
\on{SL}(2,\ZZ)$ follows from the Dedekind $\eta$-function
$$
  \eta(\tau)=q^{1/24} \prod_{n=1}^\infty (1-q^n),\quad q={e}^{2\pi i \tau}.
$$
This function is a modular form of weight $1/2$ with a multiplier
system of order $24$. Its square $\eta^2$ has weight~$1$ and is a modular
form with respect to a character $v_\eta$ of order $12$. For $g=2$ the
product
$$
  \Delta_5 (\tau) = \prod_{(m,m')\,\mathrm{even}} \Theta_{mm'}(\tau,0)
$$
of the 10 even theta characteristics is a modular form for $\on{Sp}(4,\ZZ)$
of weight 5 with respect to a character of order 2.

In \cite{GrH3} the commutator subgroups of the groups $\Gamma_{1,t}$ and
$\Gamma_{1,t}^+$ were computed. For $t\geq 1$ we put
$$
  t_1 = (t,12), \qquad t_2=(2t,12).
$$

\begin{theorem}[\cite{GrH3}]
For the commutator subgroups $\Gamma_{1,t}'$ of $\Gamma_{1,t}$ and
$(\Gamma_{1,t}^+)'$ of $\Gamma_{1,t}^+$ one obtains
\begin{enumerate}
\item[\rm (i)] $\Gamma_{1,t} / \Gamma_{1,t}' \cong \ZZ_{t_1}\times
\ZZ_{t_2}$
\item[\rm (ii)] $\Gamma_{1,t}^+ / (\Gamma_{1,t}^+)' \cong \ZZ_2\times
  \ZZ_{t_2}$.
\end{enumerate}
\end{theorem}

This was shown in \cite[Theorem 2.1]{GrH3}.

In \cite{Mu1} Mumford pointed out an interesting application of the
computation of $\on{Sp}(2,\ZZ)'$ to the Picard group of the
moduli stack $\underline{\calA}_1$. He showed that
$$
  \on{Pic}(\underline{\calA}_1)\cong \ZZ_{12}.
$$
In the same way the above theorem implies that
$$
  \on{Pic}(\underline{\calA}_2 )= \on{Pic}(\underline{\calA}_{1,1})\cong
  \ZZ \times \ZZ_2
$$
and
$$
  \on{Tors} \on{Pic}(\underline{\calA}_{1,t})=\ZZ_{t_1}\times \ZZ_{t_2}.
$$
The difference between the cases $\underline{\calA}_{1,1}$ and
$\underline{\calA}_{1,t}$, $t>1$ is that one knows that the rank of
the Picard group of $\underline{\calA}_2 = \underline{\calA}_{1,1}$ is
1, whereas the rank of the Picard group of $\underline{\calA}_{1,t}$,
$t>1$ is unknown. One only knows that it is positive. This is true for
all moduli stacks of abelian varieties of dimension $g\geq 2$, since
the bundle $L$ of modular forms of weight 1 is non-trivial. The
difference from the genus 1 case lies in the fact that there the
boundary of the Satake compactification is a divisor.

\begin{problem}
Determine the rank of the Picard group $\on{Pic}(\underline{\calA}_{1,t})$.
\end{problem}

We have already discussed Gritsenko's result which gives the existence of
weight 3 cusp forms for $\Gamma_{1,t}$ for all but finitely many values of
$t$. We call these values
$$
  t=1,2,\ldots,12,14,15,16,18,20,24,30,36
$$
the \emph{exceptional} polarizations. In many cases the results of
Gross and Po\-pescu show that weight~$3$ cusp forms indeed cannot exist.
The best possible one can hope for is the existence of weight~$3$ cusp
forms with a character of a small order. The following result is such
an existence theorem.

\begin{theorem}[\cite{GrH3}]
  Let $t$ be exceptional.
  \begin{enumerate}
  \item[\rm (i)] If $t\neq 1,2,4,5,8,16$ then there exists a weight $3$
    cusp form with respect to $\Gamma_{1,t}$ with a character of order $2$.
    \item[\rm (ii)] For $t=8,16$ there exists a weight $3$ cusp form with a
    character of order $4$.
    \item[\rm (iii)] For $t\equiv 0\mod 3, t\neq 3,9$ there exists a weight
$3$
    cusp form with a character of order $3$.
  \end{enumerate}
\end{theorem}

To every character $\chi: \Gamma_{1,t} \to \CC^\ast$ one can associate
a Siegel modular variety
$$
  \calA (\chi)=\Ker \chi \setminus \HH_2.
$$
The existence of a non-trivial cusp form of weight~$3$ with a
character $\chi$ then implies by Freitag's theorem the existence of a
differential form on a smooth projective model $\tilde\calA (\chi)$ of
$\calA(\chi)$. In particular the above result proves the existence of
abelian covers $\calA(\chi)\to\calA_{1,t}$ of small degree with
$p_g(\tilde \calA (\chi))>0$.

The proof is again an application of Gritsenko's lifting
techniques. To give the reader an idea we shall discuss the case
$t=11$ which is particularly interesting since by the result of Gross
and Popescu $\calA_{1,11}$ is unirational, but not rational. In this
case $\Gamma_{1,11}$ has exactly one character $\chi_2$. This
character has order 2. By the above theorem there is a degree 2 cover
$\calA(\chi_2)\to \calA_{1,11}$ with positive geometric genus. In
this case the lifting procedure gives us a map
$$
  \mbox{Lift:}\, J^{\rm cusp}_{3,\frac{11}{2}} ( v_\eta^{12}\times
  v_H)\to S_3 (\Gamma_{1,11}, \chi_2).
$$
Here $v_\eta$ is the multiplier system of the Dedekind $\eta$-function
and $v_\eta^{12}$ is a character of order 2. The character $v_H$ is a
character of order 2 of the integer Heisenberg group $H=H(\ZZ)$. By
$J^{\rm cusp}_{3,\frac{11}{2}} (v_\eta^2\times v_H)$ we denote the
Jacobi cusp forms of weight 3 and index $11/2$ with a character
$v_\eta^{12}\times v_H$. Similarly $S_3 (\Gamma_{1,11}, \chi_2)$ is
the space of weight 3 cusp form with respect to the group
$\Gamma_{1,11}$ and the character $\chi_2$. Recall the Jacobi theta
series
$$
  \vartheta(\tau,z)=\sum_{m\in \ZZ}\left(-\frac{4}{m}\right) q^{m^2/8}
  r^{m/2}\qquad (q=\mathrm{e}^{2\pi i \tau},\, r= \mathrm{e}^{2\pi i z})
$$
where
$$
  \left( -\frac{4}{m}\right) =\begin{cases}
    \pm 1 &\mbox{ if }m\equiv \pm 1 \mod 4\\
    0     &\mbox{ if }m\equiv 0 \mod 2.
  \end{cases}
$$
This is a Jacobi form of weight $1/2$, index $3/2$ and multiplier system
$v_\eta^3 \times v_H$. For an integer $a$ we can consider the Jacobi form
$$
  \vartheta_a (\tau,z) = \vartheta(\tau,az)\in
  J_{\frac{1}{2},\frac{1}{2}a^2}(v_\eta^3\times v_H^a).
$$
One then obtains the desired Siegel cusp form by taking
$$
F=\on{Lift}(\eta^3 \vartheta^2 \vartheta_3)\in S_3 (\Gamma_{1,11},\chi_2).
$$
Finally we want to consider the maximal abelian covering of $\calA_{1,t}$,
namely the Siegel modular threefold
$$
  \calA_{1,t}^{\rm com}= \Gamma_{1,t}' \setminus \HH_2.
$$
By $\tilde \calA_{1,t}^{\rm com}$ we denote a smooth projective model of
$\calA_{1,t}^{\rm com}$.

\begin{theorem}[\cite{GrH3}]
\begin{enumerate}
  \item[\rm (i)] The geometric genus of $\tilde\calA_{1,t}^{\rm
    com}$ is $0$ if and only if $t=1,2,4,5$.
  \item[\rm (ii)] The geometric genus of $\calA_{1,3}^{\rm com}$
    and $\calA_{1,7}^{\rm
    com}$ is $1$.
\end{enumerate}
\end{theorem}
The proof can be found as part of the proof of \cite[Theorem 3.1]{GrH3}.

At this point we should like to remark that all known construction methods
fail when one wants to construct modular forms of small weight with respect
to the groups $\Gamma_{1,t}^+$ or $\Gamma_{1,t}^\ast$. We therefore pose
the

\begin{problem}
  Construct modular forms of small weight with respect to the groups
  $\Gamma_{1,t}^+$ and $\Gamma_{1,t}^\ast$.
\end{problem}

\renewcommand{\thetheorem}{\Roman{section}.\arabic{subsection}.\arabic{theorem}}
\addtocounter{theorem}{1}
\section{Projective models}\label{IV}
In this section we describe some cases in which a Siegel
modular variety is or is closely related to an interesting projective
variety. Many of the results are very old.

\subsection{The Segre cubic}\label{IV.1}
\setcounter{theorem}{0}
Segre's cubic primal, or the {\em Segre cubic}, is the subvariety $\cS_3$ of
$\PP^5$ given by the equations
\[
\sum_{i=0}^5 x_i = \sum_{i=0}^5 x_i^3 =0
\]
in homogeneous coordinates $(x_0:\ldots:x_5)$ on $\PP^5$. Since it
lies in the hyperplane $\big(\sum x_i=0\big)\subset\PP^5$ it may be
thought of as a cubic hypersurface in $\PP^4$, but the equations as
given here have the advantage of showing that there is an action of
the symmetric group $\on{Sym}(6)$ on $\cS_3$.

These are the equations of $\cS_3$ as they are most often given in the
literature but there is another equally elegant formulation: $\cS_3$ is
given by the equations
\[
\sigma_1(x_i) = \sigma_3(x_i) =0
\]
where $\sigma_k(x_i)$ is the $k$th elementary symmetric polynomial in
the~$x_i$,
\[
\sigma_k(x_i)=\sum_{\# I=k}\prod_{i\in I}x_i.
\]
To check that these equations do indeed define $\cS_3$ it is enough to
notice that
\[
3\sigma_3(x_i)=\left(\sum x_i\right)^3-3\left(\sum x_i\right)
\left(\sum x_i^2\right)-\sum x_i^3.
\]

\begin{lemma}\label{theo23}$\cS_3$ is invariant under the action of
$\Sym(6)$ and has ten nodes, at the points equivalent to
${(1:1:1:-1:-1:-1)}$ under the $\Sym(6)$-action. This is the maximum
possible for a cubic hypersurface in $\PP^4$, and any cubic
hypersurface with ten nodes is projectively equivalent to $\cS_3$.
\end{lemma}

Many other beautiful properties of the Segre cubic and related
varieties were discovered in the nineteenth century.

The dual variety of the Segre cubic is a quartic hypersurface
$\cI_4\subset\PP^4$, the {\em Igusa quartic}. If we take homogeneous
coordinates $(y_0:\ldots:y_5)$ on $\PP^5$ then it was shown by Baker
\cite{Ba1} that $\cI_4$ is given by
\[
\sum_{i=0}^5 y_i=a^2+b^2+c^2-2(ab+bc+ca)=0
\]
where
\[
a=(y_1-y_5)(y_4-y_2),\ \  b=(y_2-y_3)(y_5-y_0)\mbox{ and
}c=(y_0-y_4)(y_3-y_1).
\]

This can also be written in terms of symmetric functions in suitable
variables as
\[
\sigma_1(x_i) = 4\sigma_4(x_i)-\sigma_2(x_i)^2 =0.
\]

This quartic is singular along ${\binom{6}{2}}=15$ lines $\ell_{ij}$,
$0\le i,j\le 5$, and ${\ell_{ij}\cap\ell_{mn}=\emptyset}$ if and only
if $\{i,j\}\cap\{m,n\}\neq\emptyset$. There are ${\frac{1}{2}}
{\binom{6}{3}}=10$ smooth quadric surfaces $Q_{ijk}$ in $\cI_4$, such
that, for instance, $\ell_{01}$, $\ell_{12}$ and $\ell_{20}$ lie in
one ruling of $Q_{012}=Q_{345}$ and $\ell_{34}$, $\ell_{45}$ and
$\ell_{53}$ lie in the other ruling. The birational map
$\cI_4\ratmap\cS_3$ given by the duality blows up the $15$ lines
$\ell_{ij}$, which resolves the singularities of~$\cI_4$, and blows
down the proper transform of each $Q_{ijk}$ (still a smooth quadric)
to give the ten nodes of~$\cS_3$.

It has long been known that if $H\subset \PP^4=\big(\sum\limits_{i=0}^5
y_i\big)$ is a hyperplane which is tangent to $\cI_4$ then
$H\cap\cI_4$ is a Kummer quartic surface. This fact provides a
connection with abelian surfaces and their moduli. The Igusa quartic
can be seen as a moduli space of Kummer surfaces. In this case,
because the polarization is principal, two abelian surfaces giving the
same Kummer surface are isomorphic and the (coarse) moduli space of
abelian surfaces is the same as the moduli space of Kummer
surfaces. This will fail in the non-principally polarized case, in
IV.3, below.

\begin{theorem}\label{theo24}$\cS_3$ is birationally equivalent to a
compactification of the moduli space $\cA_2(2)$ of principally
polarized abelian surfaces with a level-$2$ structure.
\end{theorem}

The Segre cubic is rational. An explicit birational
map $\PP^3\ratmap\cS_3$ was given by Baker~\cite{Ba1} and is presented
in more modern language in~\cite{Hun}.

\begin{cor}\label{theo25}
$\cA^*_2(2)$ is rational.
\end{cor}

A much more precise description of the relation between $\cS_3$ and
$\cA_2(2)$ is given by this theorem of Igusa.

\begin{theorem}[\cite{I2}]\label{theo26}
The Igusa compactification $\cA^*_2(2)$ of the moduli space
of principally polarized abelian surfaces with a level-$2$ structure is
isomorphic to the the blow-up $\tilde\cS_3$ of $\cS_3$ in the ten
nodes. The Satake compactification $\bar{\cA}_2(2)$
is isomorphic to~$\cI_4$, which is obtained from $\tilde\cS_3$ by
contracting $15$ rational surfaces to lines.
\end{theorem}

\begin{Proof}
The Satake compactification is $\Proj\cM\big(\Gamma_2(2)\big)$, where
$\cM(\Gamma)$ is the ring of modular forms for the group~$\Gamma$. The
ten even theta characteristics determine ten theta constants
$\theta_{m_0}(\tau),\ldots,\theta_{m_9}(\tau)$ of weight $\frac{1}{2}$
for $\Gamma_2(2)$, and $\theta^4_{m_i}(\tau)$ is a modular form of
weight~$2$ for $\Gamma_2(2)$. These modular forms determine a
map $f:\cA_2(2)\to\PP^9$ whose image actually lies in a certain
$\PP^4\subset\PP^9$. The integral closure of the subring of
$\cM\big(\Gamma_2(2)\big)$ generated by the $\theta_{m_i}^4$ is the
whole of $\cM\big(\Gamma_2(2)\big)$ and there is a quartic relation
among the $\theta_{m_i}^4$ (as well as five linear relations defining
$\PP^4\subset\PP^9$) which, with a suitable choice of basis, is the
quartic $a^2+b^2+c^2-2(ab+bc+ca)=0$. Furthermore, $f$~is an embedding
and the closure of its image is normal, so it is the Satake
compactification.\end{Proof}

The Igusa compactification is, in this context, the blow-up of the
Satake compactification along the boundary, which here consists of the
fifteen lines $\ell_{ij}$. The birational map $\cI_4\ratmap\cS_3$ does
this blow-up and also blows down the ten quadrics $Q_{ijk}$ to the ten
nodes of~$\cS_3$.

For full details of the proof see~\cite{I2}; for a more extended
sketch than we have given here and some further facts,
see~\cite{Hun}. We mention that the surfaces $Q_{ijk}$, considered as
surfaces in $\cA_2(2)$, correspond to principally polarized abelian
surfaces which are products of two elliptic curves.

Without going into details, we mention also that $\cI_4$ may be
thought of as the natural compactification of the moduli of ordered
$6$-tuples of distinct points on a conic in~$\PP^2$. Such a $6$-tuple
determines $6$~lines in $\check\PP^2$ which are all tangent to some
conic, and the Kummer surface is the double cover of $\check\PP^2$
branched along the six lines. The order gives the level-$2$ structure
(note that $\Gamma_2/\Gamma_2(2)\cong\Sp (4,\ZZ_2)\cong\Sym(6)$.) The
abelian surface is the Jacobian of the double cover of the conic
branched at the six points. On the other hand, $\cS_3$ may be thought
of as the natural compactification of the moduli of ordered $6$-tuples
of points on a line: for this, see~\cite{DO}.

The topology of the Segre cubic and related spaces has been studied by
van der Geer~\cite{vdG1} and by Lee and Weintraub~\cite{LW1},
\cite{LW2}. The method in \cite{LW1} is to show that the
isomorphism between the open parts of $\cS_3$ and $\cA_2(2)$ is
defined over a suitable number field and use the Weil conjectures.

\begin{theorem}[\cite{LW1},\cite{vdG1}]\label{theo27}
The homology of the Igusa compactification of $\cA_2(2)$ is
torsion-free. The Hodge numbers are $h^{0,0}=h^{3,3}=1$,
$h^{1,1}=h^{2,2}=16$ and $h^{p,q}=0$ otherwise.
\end{theorem}

By using the covering $\cA_2(4)\to\cA_2(2)$, Lee and
Weintraub~\cite{LW3} also prove a similar result for~$\cA_2(4)$.

\subsection{The Burkhardt quartic}\label{IV.2}
\setcounter{theorem}{0}
The {\em Burkhardt quartic} is the subvariety $\cB_4$ of $\PP^4$ given by
the equation
\[
y_0^4-y_0(y_1^3+y_2^3+y_3^3+y_4^3)+3y_1y_2y_3y_4=0.
\]
This form of degree~$4$ was found by Burkhardt~\cite{Bu} in~1888. It
is the invariant of smallest degree of a certain action of the finite
simple group $\PSp(4,\ZZ_3)$ of order~25920 on~$\PP^4$,
which arises in the study of the 27 lines on a cubic surface. In fact
this group is a subgroup of index~$2$ in the Weyl group $W(E_6)$ of~$E_6$,
which is the automorphism group of the configuration of the 27
lines. The 27 lines themselves can be recovered by solving an equation
whose Galois group is $W(E_6)$ or, after adjoining a square root of
the discriminant, $\PSp(4,\ZZ_3)$.

\begin{lemma}\label{theo28}
$\cB_4$ has forty-five nodes. Fifteen of them are equivalent to
$(1:-1:0:0:0:0)$ under the action of $\Sym(6)$ and the other
thirty are equivalent to $(1:1:\xi_3:\xi_3:\xi_3^2:\xi_3^2)$,
where $\xi_3=e^{2\pi i/3}$. This is the greatest number of nodes that
a quartic hypersurface in $\PP^4$ can have and any quartic
hypersurface in $\PP^4$ with $45$ nodes is projectively equivalent
to~$\cB_4$.
\end{lemma}

This lemma is an assemblage of results of Baker \cite{Ba2} and
de~Jong, Shepherd-Barron and Van~de~Ven \cite{JSV}: the bound on the
number of double points is the Varchenko (or spectral)
bound~\cite{Va}, which in this case is sharp.

We denote by $\theta_{\alpha\beta}(\tau)$, $\alpha,\beta\in\ZZ_3$, the
theta constants
\[
\theta_{\alpha\beta}(\tau)=\theta\left[\begin{array}{cc}0&0\\ \alpha
&\beta\end{array}\right](\tau,0)=\sum\limits_{n\in\ZZ^2}\exp\{\pi
i\,^tn\tau n+2\pi i(\alpha n_1+\beta n_2)\}
\]
where $\tau\in\HH_2$. Here we identify $\alpha\in\ZZ_3$ with
$\alpha/3\in\QQ$. The action of $\Gamma_2(1)=\Sp(4,\ZZ)$ on $\HH_2$
induces a linear action on the space spanned by these
$\theta_{\alpha\beta}$, and $\Gamma_2(3)$ acts trivially on the
corresponding projective space. Since $-1\in\Gamma_2(1)$ acts
trivially on $\HH_2$, this gives an action of
$\PSp(4,\ZZ)/\Gamma_2(3)\cong\PSp(4,\ZZ_3)$ on~$\PP^8$. The subspace
spanned by the
$y_{\alpha\beta}=\frac{1}{2}(\theta_{\alpha\beta}+\theta_{-\alpha,-\beta})$
is invariant. Burkhardt studied the ring of invariants of this
action. We put $y_0=-y_{00}$, $y_1=2y_{10}$, $y_2=2y_{01}$,
$y_3=2y_{11}$ and $y_4=2y_{1,-1}$.

\begin{theorem}[\cite{Bu},\cite{vdG2}]\label{theo29} The quartic form
$y_0^4-y_0(y_1^3+y_2^3+y_3^3+y_4^3)+3y_1y_2y_3y_4$ is an invariant, of
lowest degree, for this action. The map
\[
\tau\longmapsto(y_0:y_1:y_2:y_3:y_4)
\]
defines a map $\HH_2/\Gamma_2(3)\to\cB_4$ which extends to a birational
map {$\cA^*_2(3)\ratmap\cB_4$}.
\end{theorem}

This much is fairly easy to prove, but far more is true: van der Geer,
in \cite{vdG2}, gives a short modern proof as well as providing more
detail. The projective geometry of $\cB_4$ is better understood by
embedding it in $\PP^5$, as we did for $\cS_3$. Baker~\cite{Ba2} gives
explicit linear functions $x_0,\ldots,x_5$ of $y_0,\ldots,y_4$ such
that $\cB_4\subset\PP^5$ is given by
\[
\sigma_1(x_i) = \sigma_4(x_i) =0.
\]
The details are reproduced in~\cite{Hun}.

\begin{theorem}[\cite{To},\cite{Ba2}]\label{theo30}
$\cB_4$ is rational: consequently $\cA^*_2(3)$ is rational.
\end{theorem}

This was first proved by Todd~\cite{To}; later Baker~\cite{Ba2} gave an
explicit birational map from $\PP^3$ to~$\cB_4$.

To prove Theorem~\ref{theo29} we need to say how to recover a principally
polarized abelian surface and a level-$3$ structure from a general
point of~$\cB_4$. The linear system on a principally polarized abelian
surface given by three times the polarization is very ample, so the
theta functions $\theta_{\alpha\beta}(\tau,z)$ determine an embedding
of $A_\tau=\CC^2/\ZZ^2+\ZZ^2\tau$ ($\tau\in\HH_2$)
into~$\PP^8$. Moreover the extended Heisenberg group $G_3$ acts on the
linear space spanned by the $\theta_{\alpha\beta}$. The Heisenberg
group of level~$3$ is a central extension
\[
0\To\mu_3\To H_3\To \ZZ_3^2\To 0
\]
and $G_3$ is an extension of this by an involution~$\iota$. The
involution acts by $z\mapsto -z$ and $\ZZ_3^2$ acts by
translation by $3$-torsion points. The space spanned by
the~$y_{\alpha\beta}$ is invariant under the normalizer of the
Heisenberg group in $\on{PGL}(4,\CC)$, which is isomorphic to
$\PSp(4,\ZZ_3)$, so we get an action of this group
on $\PP^4$ and on $\cB_4\subset\PP^4$.

For a general point $p\in\cB_4$ the hyperplane in $\PP^4$ tangent to
$\cB_4$ at~$p$  meets $\cB_4$ in a quartic surface with six nodes, of
a type known as a Weddle surface. Such a surface is birational to a
unique Kummer surface (Hudson~\cite{Hud} and Jessop~\cite{Je} both
give constructions) and this is the Kummer surface of $A_\tau$.

It is not straightforward to see the level-$3$ structure in this
picture. One method is to start with a principally polarized abelian
surface $(A,\Theta)$ and embed it in $\PP^8$ by $|3\Theta|$. Then
there is a projection $\PP^8\to\PP^3$ under which the image of $A$ is
the Weddle surface, so one identifies this $\PP^3$ with the tangent
hyperplane to $\cB_4$. The Heisenberg group acts on $\PP^8$ and on
$H^0\big(\PP^8, \cO_{\PP^8}(2)\big)$, which has dimension~$45$. In
$\PP^8$, $A$ is cut out by nine quadrics in $\PP^8$. The span of these
nine quadrics is determined by five coefficients
$\alpha_0,\ldots,\alpha_4$ which satisfy a homogeneous
Heisenberg-invariant relation of degree~$4$. As the Heisenberg group
acting on $\PP^4$ has only one such relation this relation must again
be the one that defines~$\cB_4$. Thus the linear space spanned by nine
quadrics, and hence $A$ with its polarization and Heisenberg action,
are determined by a point of $\cB_4$. The fact that the two degree~$4$
relations coincide is equivalent to saying that $\cB_4$ has an unusual
projective property, namely it is self-Steinerian.

It is quite complicated to say what the level-$3$ structure means for
the Kummer surface. It is not enough to look at the Weddle surface:
one also has to consider the image of $A$ in another projection
$\PP^8\to\PP^4$, which is again a birational model of the Kummer
surface, this time as a complete intersection of type $(2,3)$ with ten
nodes. More details can be found in~\cite{Hun}.

The details of this proof were carried out by Coble~\cite{Cob}, who
also proved much more about the geometry of $\cB_4$ and the embedded
surface $A_\tau\subset\PP^8$. The next theorem is a consequence of
Coble's results.

\begin{theorem}[\cite{Cob}]\label{theo31} Let $\pi:\tilde\cB_4\to\cB_4$ be the
  blow-up of $\cB_4$ in the $45$ nodes. Then $\tilde\cB_4\cong\cA^*_2(3)$;
  the exceptional surfaces in $\tilde\cB_4$ correspond to the Humbert
  surfaces that parametrize product abelian surfaces. The Satake
  compactification is obtained by contracting the preimages of $40$ planes
  in $\cB_4$, each of which contains $9$ of the nodes.
\end{theorem}

One should compare the birational map $\cA^*_2(3)\ratmap\cB_4$ with
the birational map $\cI_4\ratmap\cS_3$ of the previous section.

By computing the zeta function of $\tilde\cB_4$ over $\FF_q$ for
$q\equiv 1$ (mod~$3$), Hoffman and Weintraub~\cite{HoW} calculated the
cohomology of~$\cA^*_2(3)$.

\begin{theorem}[\cite{HoW}]\label{theo32}
$H^i(\cA^*_2(3),\ZZ)$ is free: the odd Betti numbers are zero and
$b_2=b_4=61$.
\end{theorem}

In fact \cite{HoW} gives much more detail, describing the mixed Hodge
structures, the intersection cohomology of the Satake compactification,
the $\PSp(4,\ZZ_3)$-module structure of the cohomology and some of the
cohomology of the group $\Gamma_2(3)$. The cohomology of $\Gamma_2(3)$
was also partly computed, by another method, by MacPherson and
McConnell~\cite{McMc}, but neither result contains the other.

\subsection{The Nieto quintic}\label{IV.3}
\setcounter{theorem}{0}
The {\em Nieto quintic} $\cN_5$ is the subvariety of $\PP^5$ given in
homogeneous coordinates $x_0,\ldots,x_5$ by
\[
\sigma_1(x_i) = \sigma_5(x_i) =0.
\]
This is conveniently written as $\sum x_i =\sum \frac{1}{x_i}=0$. As
in the cases of $\cS_3$ and $\cB_4$, this form of the equation
displays the action of $\Sym(6)$ and is preferable for most purposes to
a single quintic equation in $\PP^4$. Unlike $\cS_3$ and $\cB_4$,
which were extensively studied in the nineteenth century, $\cN_5$ and
its relation to abelian surfaces was first studied only in the 1989
Ph.D. thesis of Nieto~\cite{Ni} and the paper of Barth and Nieto~\cite{BN}.

We begin with a result of van~Straten~\cite{vS}

\begin{theorem}[\cite{vS}]\label{theo33}
$\cN_5$ has ten nodes but (unlike $\cS_3$ and $\cB_4$) it also has
some non-isolated singularities. However the quintic hypersurface in
$\PP^4$ given as a subvariety of $\PP^5$ by
\[
\sigma_1 (x_i) =\sigma_5(x_i)+\sigma_2(x_i)\sigma_3(x_i) =0.
\]
has $130$ nodes and no other singularities.
\end{theorem}

This threefold and the Nieto quintic are both special elements of the
pencil
\[
\sigma_1(x_i) = \alpha\sigma_5(x_i)+\beta\sigma_2(x_i)\sigma_3(x_i) =0
\]
and the general element of this pencil has $100$ nodes. Van der
Geer~\cite{vdG2} has analysed in a similar way the pencil
\[
\sigma_1(x_i) = \alpha\sigma_4(x_i)+\beta\sigma_2(x_i)^2 =0
\]
which contains $\cB_4$ ($45$ nodes) and $\cI_4$ ($15$ singular lines)
among the special fibres, the general fibre having $30$ nodes.

No example of a quintic $3$-fold with more than $130$ nodes is known,
though the Varchenko bound in this case is $135$.

$\cN_5$, like $\cS_3$ and $\cB_4$, is related to abelian surfaces via
Kummer surfaces. The Heisenberg group $H_{2,2}$, which is a central
extension
\[
0\to\mu_2\to H_{2,2}\to \ZZ_2^4\to 0
\]
acts on $\PP^3$ via the Schr\"odinger representation on~$\CC^4$. This
is fundamental for the relation between $\cN_5$ and Kummer surfaces.

\begin{theorem}[\cite{BN}]\label{theo34}
  The space of $H_{2,2}$-invariant quartic surfaces in $\PP^3$ is
  $5$-dimensional. The subvariety of this $\PP^5$ which consists of those
  $H_{2,2}$-invariant quartic surfaces that contain a line is
  three-dimensional and its closure is projectively equivalent to~$\cN_5$.
  There is a double cover $\tilde\cN_5\to\cN_5$ such that $\tilde\cN_5$ is
  birationally equivalent to $\cA^*_{1,3}(2)$.
\end{theorem}

\begin{Proof} A general $H_{2,2}$-invariant quartic surface $X$ containing
  a line $\ell$ will contain $16$ skew lines (namely the $H_{2,2}$-orbit
  of~$\ell$). By a theorem of Nikulin~\cite{N1} this means that $X$ is the
  minimal desingularization of the Kummer surface of some abelian
  surface~$A$. The $H_{2,2}$-action on~$X$ gives rise to a level-$2$
  structure on~$A$, but the natural polarization on~$A$ is of type~$(1,3)$.
  There is a second $H_{2,2}$-orbit of lines on~$X$ and they give rise to a
  second realization of~$X$ as the desingularized Kummer surface of another
  (in general non-isomorphic) abelian surface~$\hat A$, which is in fact the
  dual of~$A$. The moduli points of $A$ and $\hat A$ (with their respective
  polarizations, but without level structures) in $\cA_{1,3}$ are related
  by $V_3(A)=\hat A$, where $V_3$ is the Gritsenko involution described in
  III.3, above.

  Conversely, given a general abelian surface~$A$ with a
  $(1,3)$-polarization and a level-$2$ structure, let $\BKm A$ be the
  desingularized Kummer surface and $\calL$ a symmetric line bundle on $A$
  in the polarization class. Then the linear system $|\calL^{\otimes 2}|^-$
  of anti-invariant sections embeds $\BKm A$ as an $H_{2,2}$-invariant
  quartic surface and the exceptional curves become lines in this
  embedding. This gives the connection between $\cN_5$ and
$\cA_{1,3}(2)$.
\end{Proof}

The double cover $\tilde\cN_5\to\cN$ is the inverse image of $\cN_5$
under the double cover of $\PP^5$ branched along the coordinate
hyperplanes.

$\cN_5$ is not very singular and therefore resembles a smooth quintic
threefold in some respects. Barth and Nieto prove much more.

\begin{theorem}[\cite{BN}]\label{theo35}
  Both $\cN_5$ and $\tilde\cN_5$ are birationally equivalent to (different)
  Calabi-Yau threefolds. In particular, the Kodaira dimension of
  $\cA^*_{1,3}(2)$ is zero.
\end{theorem}

The fundamental group of a smooth
projective model of $\cA^*_{1,3}(2)$ is isomorphic to $\ZZ_2\times\ZZ_2$
(see \cite{S1} and II.3 above).
Hence, as R.~Livn\'e has pointed out, there are four unramified covers of
such a model which
are also Calabi-Yau threefolds. In all other cases where the Kodaira
dimension of a Siegel modular variety (of dimension~$>1$) is known, the
variety is either of general type or uniruled.

It is a consequence of the above theorem that the modular
group $\Gamma_{1,3}(2)$ which defines the moduli space $\cA_{1,3}(2)$
has a unique weight-$3$ cusp form (up to a scalar). This
cusp form was determined in~\cite{GrH4}. Recall that there is a
weight-$3$ cusp form $\Delta_1$ for the group $\Gamma_{1,3}$ with
a character of order $6$. The form $\Delta_1$ has several interesting
properties, in particular it admits an infinite
product expansion and determines a generalized Lorentzian
Kac-Moody superalgebra of Borcherds type (see~\cite{GrN}).

\begin{theorem}[\cite{GrH4}]\label{theo35a}
 The modular form $\Delta_1$ is the unique weight-$3$ cusp form of
 the group $\Gamma_{1,3}(2)$.
 \end{theorem}

Using this, it is possible to give an explicit construction of a Calabi-Yau
model of $\cA_{1,3}(2)$ which does not use the projective geometry
of~\cite{BN}.

Nieto and the authors of the present survey have investigated the relation
between $\tilde\cN_5$ and $\cA^*_3(2)$ in more detail. $\cN_5$ contains
$30$ planes which fall naturally into two sets of $15$, the so-called S-
and V-planes.

\begin{theorem}[\cite{HNS1}]\label{theo36}
  The rational map $\cA^*_{1,3}(2)\ratmap\cN_5$ (which is generically
  2-to-1) contracts the locus of product surfaces to the $10$ nodes. The
  locus of bielliptic surfaces is mapped to the V-planes and the boundary
  of $\cA^*_{1,3}(2)$ is mapped to the S-planes. Thus by first blowing up
  the singular points and then contracting the surfaces in $\tilde\cN_5$
  that live over the S-planes to curves one obtains the Satake
  compactification.
\end{theorem}

In~\cite{HNS2} we gave a description of some of the degenerations that
occur over the S-planes.

One of the open problems here is to give a projective description of the
branch locus of this map. The projective geometry associated with the Nieto
quintic is much less worked out than in the classical cases of the Segre
cubic and the Burkhardt quartic.

Van Straten also calculated the Hodge numbers of a natural
desingularization of $\cN$ which is used by Barth and Nieto, by counting
points over finite fields: see~\cite{BN}.

\section{Non-principal polarizations}\label{V}

We have encountered non-principal polarizations and some of the properties
of the associated moduli spaces already. For abelian surfaces, a few of
these moduli spaces have good descriptions in terms of projective geometry,
and we will describe some of these results for abelian surfaces below. We
begin with the most famous case, historically the starting point for much
of the recent work on the whole subject.

\subsection{Type $(1,5)$ and the Horrocks-Mumford bundle}\label{V.1}
\setcounter{theorem}{0}
In this section we shall briefly describe the relation between the
Horrocks-Mumford bundle and abelian surfaces. Since this material has been
covered extensively in another survey article (see \cite {H1} and the
references quoted there) we shall be very brief here.

The existence of the Horrocks-Mumford bundle is closely related to abelian
surfaces embedded in $\PP^4$. Indeed, let $A\subset\PP^4$ be a smooth
abelian surface. Since $\omega_A={\cal O}_A$ it follows that the
determinant of the normal bundle of $A$ in $\PP^4$ is $\on{det
}N_{A/\PP^4}={\cal O}_A(5)={\cal O}_{\PP^4}(5)_{|A}$, i.e. it can be
extended to $\PP^4$. It then follows from the Serre construction (see e.g.
\cite[Theorem 5.1.1]{OSS}) that the normal bundle $N_{A/\PP^4}$ itself can
be extended to a rank $2$ bundle on $\PP^4$. On the other hand the double
point formula shows immediately that a smooth abelian surface in $\PP^4$
can only have degree $10$, so the hyperplane section is a polarization of
type $(1,5)$. Using Reider's criterion (see e.g. \cite[chapter 10, \S
4]{LB}) one can nowadays check immediately that a polarization of type
$(1,n),n\ge 5$ on an abelian surface with Picard number $\rho(A)=1$ is very
ample. The history of this subject is, however, quite intricate. Comessatti
proved in 1916 that certain abelian surfaces could be embedded in $\PP^4$.
He considered a 2-dimensional family of abelian surfaces, namely those
which have real multiplication in $\QQ(\sqrt{5})$. His main tool was theta
functions. His paper~\cite{Com} was later forgotten outside the Italian
school of algebraic geometers. A modern account of Comessatti's results
using, how\-ever, a different language and modern methods was later given
by Lange \cite{L} in 1986. Before that Ramanan \cite{R} had proved a
criterion for a $(1,n)$-polarization to be very ample.  This criterion
applies to all $(1,n)$-polarized abelian surfaces $(A,H)$ which are cyclic
$n$-fold covers of a Jacobian. In particular this also gives the existence
of abelian surfaces in $\PP^4$. The remaining cases not covered by
Ramanan's paper were treated in~\cite{HL}.

With the exception of Comessatti's essentially forgotten paper, none of
this was available when Horrocks and Mumford investigated the existence of
indecomposable rank 2 bundles on $\PP^4$. Although they also convinced
themselves of the existence of smooth abelian surfaces in $\PP^4$ they then
presented a construction of their bundle $F$ in \cite{HM} in cohomological
terms, i.e. they constructed $F$ by means of a {\em monad}. A monad is a
complex
$$
(M)\qquad A\stackrel{p}{\longrightarrow} B \stackrel{q}{\longrightarrow}C
$$
where $A,B$ and $C$ are vector bundles, $p$ is injective as a map of
vector bundles, $q$ is surjective and $q\circ p=0$. The cohomology of $(M)$
is
$$
F=\Ker q/\Im p
$$
which is clearly a vector bundle. The Horrocks-Mumford bundle can be
given by a monad of the form
$$
V\otimes {\cal O}_{\PP^4}(2)\stackrel{p}{\rightarrow}
2\bigwedge\nolimits^2T_{\PP^4}\stackrel{q}{\rightarrow}V^{\ast}\otimes{\cal
O}_{\PP^4}(3)
$$
where $V=\CC^5$ and $\PP^4=\PP(V)$. The difficulty is to write down the
maps $p$ and $q$. The crucial ingredient here is the maps
$$
\begin{array}{lcllcl}
f^+&:& V\longrightarrow\bigwedge\nolimits^2 V,&f^+(\sum v_i e_i)&=&
\sum v_i e_{i+2}\wedge e_{i+3}\\
f^{-1}&:& V\longrightarrow\bigwedge\nolimits^2 V,&f^-(\sum v_i e_i)&=&
\sum v_i e_{i+1}\wedge e_{i+4}
\end{array}
$$
where $(e_i)_{i\in\ZZ_5}$ is the standard basis of $V=\CC^5$ and indices
have to be read cyclically. The second ingredient is the Koszul complex on
$\PP^4$, especially its middle part
$$
\diagram
\bigwedge^2 V\otimes{\cal O}_{\PP^4}(1)\rrto^{\wedge s}\drto^{p_{0}} & &
\bigwedge^3 V\otimes{\cal O}_{\PP^4}(2)\\
& \bigwedge^2 T_{\PP^4}(-1)\urto^{q_{0}}
\enddiagram
$$
where $s:{\cal O}_{\PP^4}(-1)\rightarrow V\otimes{\cal O}_{\PP^4}$ is
the tautological bundle map. The maps $p$ and $q$ are then given by
$$
\begin{array}{lclcl}
p:V\otimes{\cal O}_{\PP^4}(2) &\stackrel{(f_+, f_-)}{\longrightarrow}&
2\bigwedge^2V\otimes {\cal O}_{\PP^4}(2)&\stackrel{2
p_{0}(1)}{\longrightarrow}&2\bigwedge^2 T_{\PP^4}\\
q:2\bigwedge^2 T_{\PP^4} &\stackrel{2q_{0}(1)}{\longrightarrow}&
2\bigwedge^3V\otimes {\cal
O}_{\PP^4}(3)&\stackrel{-(^tf_-, ^tf_+)}{\longrightarrow}&
V^{\ast}\otimes{\cal O}_{\PP^4}(3).
\end{array}
$$
Once one has come up with these maps it is not difficult to check that
$p$ and $q$ define a monad. Clearly the cohomology $F$ of this monad
is a rank $2$ bundle and it is straightforward to calculate its Chern
classes to be
$$
c (F)=1+5h+10h^2
$$
where $h$ denotes the hyperplane section. Since this polynomial is
irreducible over the integers it follows that $F$ is indecomposable.

One of the remarkable features of the bundle $F$ is its symmetry
group. The {\em Heisenberg group} of level $n$ is the
subgroup $H_n$ of $\SL(n,\CC)$ generated by the automorphisms
$$
\sigma:\quad e_i\mapsto e_{i-1},\quad \tau:e_i\mapsto \varepsilon^i e_i\quad
(\varepsilon=e^{2\pi i/n}).
$$
Since $[\sigma, \tau]=\varepsilon\cdot\mbox{id}_V$ the group $H_n$ is a
central extension
$$
0\rightarrow \mu_n\rightarrow H_n\rightarrow \ZZ_n\times \ZZ_n\rightarrow 0.
$$
Let $N_5$ be the normalizer of the Heisenberg group $H_5$ in
$\SL(5,\CC)$. Then $N_5/H_5\cong \SL(2,\ZZ_5)$ and $N_5$ is in fact a
semi-direct product
$$
N_5\cong H_5\rtimes \SL(2,\ZZ_5).
$$
Its order is $|N_5|=|H_5|\cdot|\SL(2,\ZZ_5)|=125\cdot 120=15,000$. One
can show that $N_5$ acts on the bundle $F$ and that it is indeed its
full symmetry group~\cite{De}.

The Horrocks-Mumford bundle is {\em stable}. This follows since
$F(-1)=F\otimes {\cal O}_{\PP^4}(-1)$ has $c_1(F(-1))=3$ and
$h^{0}(F(-1))=0$. Indeed $F$ is the unique stable rank $2$ bundle with
$c_1=5$ and $c_2=10$ \cite{DS}. The connection with abelian surfaces
is given via sections of $F$. Since $F(-1)$ has no sections every
section $0\neq s\in H^{0}(F)$ vanishes on a surface whose degree is
$c_2(F)=10$.

\begin{proposition}\label{theo37}
For a general section $s\in H^{0}(F)$ the zero-set $X_s=\{s=0\}$ is a
smooth abelian surface of degree $10$.
\end{proposition}

\begin{Proof} \cite[Theorem 5.1]{HM}. The crucial point is to prove that
$X_s$ is smooth. The vector bundle $F$ is globally generated outside
$25$ lines $L_{ij}$ in $\PP^4$. It therefore follows from Bertini that
$X_s$ is smooth outside these lines. A calculation in local
coordinates then shows that for general $s$ the surface $X_s$ is also
smooth where it meets the lines $L_{ij}$. It is then an easy
consequence of surface classification to show that $X_s$ is
abelian.
\end{Proof}

In order to establish the connection with moduli spaces it is useful
to study the space of sections $H^{0}(F)$ as an $N_5$-module. One can
show that this space is $4$-dimensional and that the Heisenberg group
$H_5$ acts trivially on $H^{0}(F)$. Hence $H^{0}(F)$ is an
$\SL(2,\ZZ_5)$-module. It turns out that the action of
$\SL(2,\ZZ_5)$ on $H^{0}(F)$ factors through an action of
$\PSL(2,\ZZ_5)\cong A_5$ and that as an $A_5$-module $H^{0}(F)$
is irreducible. Let $U\subset
\PP^3=\PP(H^{0}(F))$ be the open set parametrising {\em smooth}
Horrocks-Mumford surfaces $X_s$. Then $X_s$ is an abelian surface
which is fixed under the Heisenberg group $H_5$. The action of $H_5$
on $X_s$ defines a canonical level-$5$ structure on $X_5$. Let ${\cal
A}_{1,5}^{\lev}$ be the moduli space of triples
$(A,H,\alpha)$ where $(A,H)$ is a $(1,5)$-polarized abelian surface
and $\alpha$ a canonical level structure and denote by
$\opn{\cal A}_{1,5}^{\lev}$ the open part where
the polarization $H$ is very ample. Then the above discussion leads to

\begin{theorem}[\cite{HM}]\label{theo38}
The map which associates to a section $s$ the Hor\-rocks-Mumford surface
$X_s=\{s=0\}$ induces an isomorphism of $U$ with $\opn{\cal
A}_{1,5}^{\lev}$. Under this isomorphism the action of
$\PSL(2,\ZZ_5)=A_5$ on $U$ is identified with the action of
$\PSL(2,\ZZ_5)$ on ${\cal A}_{1,5}^{\lev}$ which
permutes the canonical level structures on a $(1,5)$-polarized
abelian surface. In particular ${\cal A}_{1,5}^{\lev}$
is a rational variety.
\end{theorem}

\begin{Proof} \cite[Theorem 5.2]{HM}.
\end{Proof}

The inverse morphism
$$
\varphi:{\opn{\cal A}}_{1,5}^{\lev}\rightarrow U\subset
\PP(H^{0}(F))=\PP^3
$$
can be extended to a morphism
$$
{\tilde{\varphi}}:(\cA_{1,5}^\lev)^*\rightarrow \PP(H^{0}(F))
$$
where $(\cA_{1,5}^\lev)^*$ denotes the Igusa (=Voronoi)
compactification of ${\cal A}_{1,5}^{\lev}$. This extension can also
be understood in terms of degenerations of abelian surfaces. Details
can be found in \cite{HKW2}.

\subsection{Type $(1,7)$}\label{V.2}
\setcounter{theorem}{0}
The case of type $(1,7)$ was studied by Manolache and Schreyer
\cite{MS} in 1993. We are grateful to them for making some
private notes and a draft version of \cite{MS} available to us and
answering our questions. Some of their results have also been found by Gross
and Popescu~\cite{GP1}, \cite{GP2} and by Ranestad: see also~\cite{S-BT}.

\begin{theorem}[\cite{MS}]\label{theo39}
$\cA^{\lev}_{1,7}$ is rational, because it is birationally equivalent
to a Fano variety of type~$V_{22}$.
\end{theorem}

\begin{Proof} We can give only a sketch of the proof here.
For a general abelian surface $A$ with a polarization of type $(1,7)$
the polarization is very ample and embeds $A$ in $\PP^6$. In the
presence of a canonical level structure the $\PP^6$ may be thought of
as $\PP(V)$ where $V$ is the Schr\"odinger representation of the
Heisenberg group~$H_7$. We also introduce, for $j\in\ZZ_7$, the
representation $V_j$, which is the Schr\"odinger representation
composed with the automorphism of $H_7$ given by $e^{2\pi i/7}\mapsto
e^{6\pi ij/7}$. These can also be thought of as representations of the
extended Heisenberg group $G_7$, the extension of $H_7$ by an extra
involution coming from $-1$ on~$A$. The representation $S$ of $G_7$ is
the character given by this involution (so $S$ is trivial on $H_7$).

It is easy to see that $A\subset\PP^6$ is not
contained in any quadric, that is $H^0\big(\cI_A(2)\big)=0$, and from
this it follows that there is an $H_7$-invariant resolution
\begin{eqnarray*}
0\gets\cI_A\gets 3V_4\otimes\cO(-3)\gets 7V_1\otimes\cO(-4)\gets
6V_2\otimes\cO(-5)\\
\mbox{               }\gets 2V\otimes\cO(-6)\oplus\cO(-7)\gets 2\cO(-7)\gets 0.
\end{eqnarray*}
By using this and the Koszul complex one obtains a symmetric
resolution
\begin{eqnarray*}
0\gets\cO_A\gets \cO{\buildrel\beta \over \gets}
3V_4\otimes\cO(-3){\buildrel \alpha \over \gets} 2S\otimes\Omega^3
{\buildrel{\alpha'}\over\gets}
3V_1\otimes\cO(-4){\buildrel{\beta'}\over\gets} \cO(-7)\gets 0.
\end{eqnarray*}
This resolution is $G_7$-invariant. Because of the $G_7$-symmetry,
$\alpha$ can be described by a $3\times 2$~matrix $X$ whose entries
lie in a certain $4$-dimensional space~$U$, which is a module for
$\SL(2,\ZZ_7)$. The symmetry of the resolution above amounts to saying
that $\alpha'$ is given by the matrix $X'=\left(\begin{array}{cc}
    0&1\\ -1&0 \end{array}\right){}^tX$, and the complex tells us that
$\alpha\alpha'=0$. The three $2\times 2$ minors of $X$ cut out a
twisted cubic curve $C_A$ in $\PP(U^\vee)$ and because of the
conditions on $\alpha$ the ideal $I_A$ of this cubic is annihilated by
the differential operators
\begin{eqnarray*}
\Delta_1&=&\frac{\partial^2}{\partial u_0\partial u_1}
-{\textstyle\frac{1}{2}}\frac{\partial^2}{\partial u_2^2},\\
\Delta_2&=&\frac{\partial^2}{\partial u_0\partial u_2}
-{\textstyle\frac{1}{2}}\frac{\partial^2}{\partial u_3^2},\\
\Delta_3&=&\frac{\partial^2}{\partial u_0\partial u_3}
-{\textstyle\frac{1}{2}}\frac{\partial^2}{\partial u_1^2}
\end{eqnarray*}
where the $u_i$ are coordinates on~$U$.

This enables one to recover the abelian surface $A$ from~$C_A$. If we write
$R=\CC[u_0,u_1,u_2,u_3]$ then we have a complex (the Hilbert-Burch complex)
$$
0\longleftarrow R/I_A\longleftarrow R\longleftarrow
R(-2)^{\oplus 3}{\buildrel X \over \longleftarrow}
R(-3)^{\oplus 2}\longleftarrow 0.
$$
It is exact, because otherwise one can easily calculate the syzygies of
$I_A$ and see that they cannot be the syzygies of any ideal annihilated by
the three~$\Delta_i$. So $I_A$ determines $\alpha$ (up to conjugation) and
the symmetric resolution of~$\cO_A$ can be reconstructed from~$\alpha$.

Let $H_1$ be the component of the Hilbert scheme parametrising twisted
cubic curves. For a general net of quadrics $\delta\subset\check\PP(U)$ the
subspace $H(\delta)\subset H_1$ consisting of those cubics annihilated by
$\delta$ is, by a result of
Mukai~\cite{Muk}, a smooth rational Fano $3$-fold of genus~$12$, of the
type known as $V_{22}$. To check that this is so in a particular case it is
enough to show that $H(\delta)$ is smooth. We must do so for
$\delta=\Delta=\on{Span}(\Delta_1,\Delta_2\Delta_3)$. Manolache and
Schreyer show that $H(\Delta)$ is isomorphic to the
space $\on{VSP}\big(\bar{X}(7),6\big)$ of polar hexagons to the Klein
quartic curve (the modular curve~$\bar{X}(7)$):
$$
\on{VSP}\big(\bar{X}(7),6\big)=\big\{\{l_1,\ldots,l_6\}
\subset\on{Hilb}^6(\check\PP^2)\mid \sum
l_i^4=x_0^3x_1+x_1^3x_2+x_2^3x_0\big\}.
$$
(To be precise we first consider all $6$-tuples $(l_1,\ldots,l_6)$ where
the $l_i$ are pairwise different with the above property and then take the
Zariski-closure in the Hilbert scheme.)
It is known that $\on{VSP}\big(\bar{X}(7),6\big)$ is smooth, so we are done.
\end{Proof}

Manolache and Schreyer also give an explicit rational parametrization of
$\on{VSP}\big(\bar{X}(7),6\big)$ by writing down equations for the abelian
surfaces. They make the interesting observation that this rational
parametrization is actually defined over the rational numbers.

\subsection{Type $(1,11)$}\label{V.3}
\setcounter{theorem}{0}
The spaces $\cA_{1,d}^{\lev}$ for small~$d$ are studied by
Gross and Popescu, \cite{GP1}, \cite{GP1a}, \cite{GP2}, \cite{GP3}.
In particular,
in~\cite{GP1a}, they obtain a description of~$\cA_{1,11}^{\lev}$.

\begin{theorem}[\cite{GP1a}]\label{theo40}
There is a rational map
$
\Theta_{11}:\cA_{1,t}^{\lev}\ratmap\Gr(2,6)
$
which is birational onto its image. The closure
of $\Im\Theta_{11}$ is a smooth linear section of $\Gr(2,6)$ in the
Pl\"ucker embedding and is birational to the Klein cubic in
$\PP^4$. In particular $\cA_{1,11}^{\lev}$ is unirational but
not rational.
\end{theorem}

The Klein cubic is the cubic hypersurface in $\PP^4$ with the equation
\[
\sum_{i=0}^4 x_i^2x_{i+1}=0
\]
with homogeneous coordinates $x_i$, $i\in\ZZ_5$. It is smooth, and all
smooth cubic hypersurfaces are unirational but not
rational~\cite{CG},~\cite{IM}.

The rational map $\Theta_{11}$ arises in the following way. For a
general abelian surface $A$ in $\cA_{1,11}^{\lev}$, the
polarization (which is very ample) and the level structure determine
an $H_{11}$-invariant embedding of $A$ into $\PP^{10}$. The action of
$-1$ on $A$ lifts to $\PP^{10}=\PP (H^0(\calL))$ and the $(-1)$-eigenspace
of this action on $H^0(\calL)$ (where $\calL$ is a symmetric bundle in the
polarizing class) determines a $\PP^4$, called
$\PP^-\subset\PP^{10}$. We choose coordinates $x_0,\ldots,x_{10}$ on
$\PP^{10}$ with indices in $\ZZ_{11}$ such that $x_1,\ldots,x_5$ are
coordinates on $\PP^-$, so that on $\PP^-$ we have $x_0=0$,
$x_i=-x_{-i}$. The matrix $T$ is defined to be the restriction of $R$
to $\PP^-$, where
\[
R_{ij}=x_{j+i}x_{j-i}, \qquad 0\le i,j\le 5
\]
(This is part of a larger matrix which describes the action on
$H^0\big(\cO_{\PP^{10}}(2)\big)$ of $H_{11}$.) The matrix $T$ is
skew-symmetric and non-degenerate at a general point of
$\PP^-$. However, it turns out that for a general
$A\in\cA_{1,11}^{\lev}$ the rank of $T$ at a general point $x\in
A\cap\PP^-$ is~$4$. For a fixed $A$, the kernel of $T$ is independent
of the choice of $x$ (except where the dimension of the kernel jumps),
and this kernel is the point $\Theta_{11}(A)\in\Gr(2,6)$.

>From the explicit matrix~$R$, finally, Gross and Popescu obtain the
description of the closure of $\Im\Theta_{11}$ as being the
intersection of $\Gr(2,6)$ with five hyperplanes in Pl\"ucker
coordinates. The equation of the Klein cubic emerges directly (as a
$6\times 6$ Pfaffian), but it is a theorem of Adler~\cite{AR} that the
Klein cubic is the only degree~$3$ invariant of $\PSL(2,\ZZ_{11})$
in~$\PP^4$.

\subsection{Other type $(1,t)$ cases}\label{V.4}
\setcounter{theorem}{0}
The results of Gross and Popescu for $t=11$ described above are part
of their more general results about $\cA_{1,t}^{\lev}$ and $\cA_{1,t}$
for $t\ge 5$.
In the series of papers \cite{GP1}--\cite{GP3} they prove the following
(already stated above as Theorem~\ref{theo18}).

\begin{theorem}[\cite{GP1},\cite{GP1a},\cite{GP2},\cite{GP3}]\label{theo41}
${\cal A}_{1,t}^{\lev}$ is rational for $6\le
t\le 10$ and $t=12$ and unirational, but not rational, for $t=11$.
Moreover the variety ${\cal A}_{1,t}$ is unirational for $t=14, 16,
18$ and $20$.
\end{theorem}

The cases have a different flavour depending on whether $t$ is even or
odd. For odd $t=2d+1$ the situation is essentially as described for
$t=11$ above: there is a rational map
$\Theta_{2d+1}:\cA_{1,t}^{\lev}\ratmap\Gr(d-3,d+1)$, which can be
described in terms of matrices or by saying that $A$ maps to the
$H_t$-subrepresentation $H^0\big(\cI_A(2)\big)$ of
$H^0\big(\cO_A(2)\big)$. In other words, one embeds $A$ in $\PP^{t-1}$
and selects the $H_t$-space of quadrics vanishing along~$A$.

\begin{theorem}[\cite{GP1}]\label{theo42}
If $t=2d+1\ge 11$ is odd then the homogeneous ideal of a general
$H_t$-invariant abelian surface in $\PP^{t-1}$ is generated by
quadrics; consequently $\Theta_{2d+1}$ is birational onto its image.
\end{theorem}

For $t=7$ and $t=9$ this is not true: however, a detailed analysis is
still possible and is carried out in~\cite{GP2} for $t=7$ and
in~\cite{GP1a} for $t=9$. For $t\ge 13$ it is a
good description of the image of $\Theta_t$ that is lacking. Even for
$t=13$ the moduli space is not unirational and for large $t$ it is of
general type (at least for $t$ prime or a prime square).

For even $t=2d$ the surface $A\subset\PP^{t-1}$ meets
$\PP^-=\PP^{d-2}$ in four distinct points (this is true even for many
degenerate abelian surfaces). Because of the $H_t$-invariance these
points form a $\ZZ_2\times\ZZ_2$-orbit and there is therefore a
rational map $\Theta_{2d}:\cA_{1,t}^{\lev}\ratmap\PP^-/(\ZZ_2\times\ZZ_2)$.

\begin{theorem}[\cite{GP1}]\label{theo43}
If $t=2d\ge 10$ is even then the homogeneous ideal of a general
$H_t$-invariant abelian surface in $\PP^{t-1}$ is generated by
quad\-rics (certain Pfaffians) and $\Theta_{2d}$ is birational onto its
image.
\end{theorem}

To deduce Theorem \ref{theo41} from Theorem \ref{theo42} and
Theorem \ref{theo43} a careful analysis of each case is necessary:
for $t=6,8$ it is again the case that $A$ is not cut out by quadrics
in $\PP^{t-1}$. In those cases when rationality or unirationality can
be proved, the point is often that there are pencils of abelian
surfaces in suitable Calabi-Yau $3$-folds and these give rise to
rational curves in the moduli spaces. Gross and Popescu use these
methods in~\cite{GP1a} ($t=9,11$), \cite{GP2} ($t=6,7,8$ and $10$),
and \cite{GP3} ($t=12$) to obtain detailed
information about the moduli spaces $\cA_{1,t}^{\lev}$.
In~\cite{GP3} they also consider the spaces $\cA_{1,t}$
for $t=14,16,18$ and $20$.

\section{Degenerations}\label{VI}

The procedure of toroidal compactification described in \cite{AMRT}
involves making many choices. Occasionally there is an obvious
choice. For moduli of abelian surfaces this is usually the case, or
nearly so, since one has the Igusa compactification (which is the
blow-up of the Satake compactification along the boundary) and all
known cone decompositions essentially agree with this one. But
generally toroidal compactifications are not so simple. One has to
make further modifications in order to obtain acceptably mild
singularities at the boundary. Ideally one would like to do this in a
way which is meaningful for moduli, so as to obtain a space which
represents a functor described in terms of abelian varieties and
well-understood degenerations. The model, of course, is the
Deligne-Mumford compactification of the moduli space of curves.

\subsection{Local degenerations}\label{VI.1}
\setcounter{theorem}{0}
The first systematic approach to the local problem of constructing
degenerations of polarized abelian varieties is Mumford's paper
\cite{Mu2} (conveniently reprinted as an appendix to
\cite{FC}). Mumford specifies degeneration data which determine a
family $G$ of semi-abelian varieties over the spectrum $S$ of a
complete normal ring $R$. Faltings and Chai \cite{FC} generalized this
and also showed how to recover the degeneration data from such a
family. This semi-abelian family can then be compactified: in fact,
Mumford's construction actually produced the compactification first
and the semi-abelian family as a subscheme. However, although $G$ is
uniquely determined, the compactification is non-canonical. We may as
well assume that $R$ is a DVR and that $G_\eta$, the generic fibre, is
an abelian scheme: the compactification then amounts to compactifying
the central fibre $G_0$ in some way.

Namikawa (see for instance \cite{Nam3} for a concise account) and
Nakamura \cite{Nak1} used toroidal methods to construct
natural compactifications in the complex-analytic category, together
with proper degenerating families of so-called stable quasi-abelian
varieties. Various difficulties, including non-reduced fibres,
remained, but more recently Alexeev and Nakamura
\cite{Ale2}, \cite{AN},
have produced a more satisfactory and simpler theory. We describe
their results below, beginning with their simplified version of the
constructions of Mumford and of Faltings and Chai. See \cite{FC},
\cite{Mu2} or \cite{AN} itself for more.

$R$ is a complete DVR with maximal ideal~$I$, residue field $k=R/I$
and field of fractions~$K$. We take a split torus $\tilde G$ over $R$
with character group $X$ and let $\tilde G(K)\cong(K^*)^g$ be the
group of $K$-valued points of~$\tilde G$. A set of periods is simply a
subgroup $Y\subset\tilde G(K)$ which is isomorphic to~$\ZZ^g$. One can
define a polarization to be an injective map $\phi:Y\to X$ with
suitable properties.

\begin{theorem}[\cite{Mu2},\cite{FC}]\label{theo46} There is a quotient
$G=\tilde G/Y$ which is a semi-abelian scheme over~$S$: the generic
fibre $G_\eta$ is an abelian scheme over $\Spec K$ with a polarization
(given by a line bundle $\calL_\eta$ induced by~$\phi$).
\end{theorem}

This is the special case of maximal degeneration, when $G_0$ is a
torus over~$k$. In practice one starts not with $\tilde G$ but with the
generic fibre $G_\eta$. According to the semistable reduction theorem
there is always a semi-abelian family $G\to S$ with generic fibre
$G_\eta$, but in order to construct a uniformization $G=\tilde G/Y$ as
above we have to allow $\tilde G$ to have an abelian part. Then
$\tilde G$ and $G_0$, instead of being split tori, are Raynaud
extensions, that is, extensions of tori by abelian schemes, over $R$
and $k$ respectively. The extra work entailed by this is carried out
in \cite{FC} but the results, though a little more complicated to
state, are essentially the same as in the case of maximal degeneration

Mumford's proof also provides a projective degeneration, in fact a
wide choice of projective degenerations, each containing $G$ as an
open subscheme.

\begin{theorem}[\cite{Mu2},\cite{Ch},\cite{FC},\cite{AN}]\label{theo47}
There is an integral scheme $\tilde P$,
locally of finite type over $R$, containing $\tilde G$ as an open
subscheme, with an ample line bundle $\tilde\calL$ and an action of $Y$
on $(\tilde P, \tilde\calL)$. There is an $S$-scheme $P=\tilde P/Y$,
projective over~$S$, with $P_\eta\cong G_\eta$ as polarized varieties,
and $G$ can be identified with an open subscheme of~$P$.
\end{theorem}

Many technical details have been omitted here. $\tilde P$ has to
satisfy certain compatibility and completeness conditions: of these,
the most complicated is a completeness condition which is used in
\cite{FC} to prove that each component of the central fibre $P_0$ is
proper over~$k$. Alexeev and Nakamura make a special choice of $\tilde
P$ which, among other merits, enables them to dispense with this
condition because the properness is automatic.

The proof of \ref{theo46}, in the version given by Chai \cite{Ch} involves
implicitly) writing down theta functions on $\tilde G(K)$ in order to check
that the generic fibre is the abelian scheme $G_\eta$. These theta
functions can be written (analogously with the complex-analytic case) as
Fourier power series convergent in the $I$-adic topology, by taking
coordinates $w_1,\ldots,w_g$ on $\tilde G(K)$ and setting
$$
\theta=\sum_{x \in X}\sigma_x(\theta)w^x
$$
with $\sigma_x(\theta)\in K$. In particular theta functions
representing elements of $H^0(G_\eta, \calL_\eta)$ can be written this
way and the coefficients obey the transformation formula
$$
\sigma_{x+\phi(y)}(\theta)=a(y)b(y,x)\sigma_x(\theta)
$$
for suitable functions $a:Y\to K^*$ and $b:Y\times X\to K^*$.

For simplicity we shall assume for the moment that the polarization is
principal: this allows us to identify $Y$ with $X$ via~$\phi$ and also
means that there is only one theta function,~$\vartheta$. The general
case is only slightly more complicated.

These power series have $K$ coefficients and converge in the $I$-adic
topology but their behaviour is entirely analogous to the familiar
complex-analytic theta functions. Thus there are cocycle conditions on $a$
and $b$ and it turns out that $b$ is a symmetric bilinear form on $X\times
X$ and $a$ is an inhomogeneous quadratic form. Composing $a$ and $b$ with
the valuation yields functions $A:X\to\ZZ$, $B:X\times X\to\ZZ$, and they
are related by
$$
A(x)=\frac{1}{2}B(x,x)+\frac{rx}{2}
$$
for some $r\in\NN$. We fix a parameter $s\in R$, so $I=sR$.

\begin{theorem}[\cite{AN}]\label{theo48} The normalization of $\Proj
R[s^{A(x)}w^x\theta; x\in X]$ is a relatively complete model $\tilde
P$ for the maximal degeneration of principally polarized abelian
varieties associated with~$G_\eta$.
\end{theorem}

Similar results hold in general. The definition of $\tilde P$ has to
be modified slightly if $G_0$ has an abelian part. If the polarization
is non-principal it may be necessary to make a ramified base change
first, since otherwise there may not be a suitable extension of
$A:Y\to\ZZ$ to $A:X\to\ZZ$. Even for principal polarization it may be
necessary to make a base change if we want the central fibre $P_0$ to
have no non-reduced components.

The proof of Theorem \ref{theo48} depends on the observation that the ring
$$
R[s^{A(x)}w^x\theta; x\in X]
$$
is generated by monomials. Consequently
$\tilde P$ can be described in terms of toric geometry. The quadratic
form $B$ defines a {\em Delaunay decomposition} of $X\otimes\RR=X_\RR$. One
of the many ways of describing this is to consider the paraboloid in
$\RR e_0\oplus X_\RR$ given by
$$
x_0=a(x)=\frac{1}{2} B(x,x)+\frac{rx}{2},
$$
and the lattice $M=\ZZ e_0\oplus X$. The convex hull of the points of
the paraboloid with $x\in X$ consists of countably many facets and
the projections of these facets on $X_\RR$ form the Delaunay
decomposition. This decomposition determines~$\tilde P$. It is
convenient to express this in terms of the {\em Voronoi decomposition}
${\mathrm{Vor}}_B$ of $X_\RR$ which is dual to the Delaunay
decomposition in the sense that there is a 1-to-1 inclusion-reversing
correspondence between (closed) Delaunay and Voronoi cells. We
introduce the map $dA:X_\RR\to X_\RR^*$ given by
$$
dA(\xi)(x)=B(\xi, x)+\frac{rx}{2}.
$$

\begin{theorem}[\cite{AN}]\label{theo49} $\tilde P$ is the torus embedding
over~$R$ given by the lattice $N=M^*\subset \RR e_0^*\oplus X^*_\RR$
and the fan $\Delta$ consisting of $\{ 0\}$ and the cones on the
polyhedral cells making up $\big(1,-dA({\mathrm{Vor}}_B)\big)$.
\end{theorem}

Using this description, Alexeev and Nakamura check the required
properties of $\tilde P$ and prove Theorem \ref{theo48}. They also obtain a
precise description of the central fibres $\tilde P_0$ (which has no
non-reduced components if we have made a suitable base change)
and~$P_0$ (which is projective). The polarized fibres $(p_0,\calL_0)$
that arise are called {\em stable quasi-abelian varieties}, as in
\cite{Nak1}. In the principally polarized case $P_0$ comes with a
Cartier divisor $\Theta_0$ and $(P_0,\Theta_0)$ is called a {\em stable
quasi-abelian pair}. We refer to \cite{AN} for a precise intrinsic
definition, which does not depend on first knowing a degeneration that
gives rise to the stable quasi-abelian variety. For our purposes all
that matters is that such a characterization exists.

\subsection{Global degenerations and compactification}\label{VI.2}
\setcounter{theorem}{0}
Alexeev, in \cite{Ale2}, uses the infinitesimal degenerations that we
have just been considering to tackle the problem of canonical global
moduli. For simplicity we shall describe results of \cite{Ale2} only in the
principally polarized case.

We define a {\em semi-abelic variety\/} to be a normal variety $P$
with an action of a semi-abelian variety $G$ having only finitely many
orbits, such that the stabilizer of the generic point of $P$ is a
connected reduced subgroup of the torus part of $G$. If $G=A$ is
actually an abelian variety then Alexeev refers to $P$ as an {\em
abelic variety\/}: this is the same thing as a torsor for the abelian
variety $A$. If we relax the conditions by allowing $P$ to be
semi-normal then $P$ is called a {\em stable semi-abelic variety\/} or SSAV.

A {\em stable semi-abelic pair\/} $(P,\Theta)$ is a projective SSAV together
with an effective ample Cartier divisor $\Theta$ on $P$ such that $\Theta$
does not contain any $G$-orbit. The degree of the corresponding
polarization is $h^0(\cO_P(\Theta))$, and $P$ is said to be principally
polarized if the degree of the polarization is~$1$. If $P$ is an
abelic variety then $(P,\Theta)$ is called an {\em abelic pair}.

\begin{theorem}[\cite{Ale2}]\label{theo50}
The categories ${\underline A}_g$ of $g$-dimensional principally
polarized abelian varieties and ${\underline{AP}}_g$ of principally
polarized abelic pairs are naturally equivalent. The
corresponding coarse moduli spaces $\cA_g$ and $\cA\cP_g$ exist as
separated schemes and are naturally isomorphic to each other.
\end{theorem}

Because of this we may as well compactify ${\cA\cP}_g$ instead of
$\cA_g$ if that is easier. Alexeev carries out this program in
\cite{Ale2}. In this way, he obtains a proper algebraic space
$\overline{\cA\cP}_g$ which is a coarse moduli space for
stable semi-abelic pairs.

\begin{theorem}[\cite{Ale2}]\label{theo51}
  The main irreducible component of $\overline{\cA\cP}_g$ (the
  component that contains $\cA\cP_g=A_g$) is isomorphic to the Voronoi
  compactification $\cA_g^*$ of~$\cA_g$. Moreover, the
  Voronoi compactification in this case is projective.
\end{theorem}

The first part of Theorem~\ref{theo51} results from a careful comparison of
the respective moduli stacks. The projectivity, however, is proved by
elementary toric methods which, in view of the results of \cite{FC},
work over $\Spec\ZZ$.

In general $\overline{\cA\cP}_g$ has other components, possibly of
very large dimension. Alexeev has examined these components and the
SSAVs that they parametrize in~\cite{Ale3}

Namikawa, in \cite{Nam1}, already showed how to attach a stable
quasi-abelian variety to a point of the Voronoi compactification.
Namikawa's families, however, have non-reduced fibres and require the
presence of a level structure: a minor technical alteration (a base
change and normalization) has to be made before the construction works
satisfactorily. See \cite{AN} for this
and also for an alternative construction using explicit local families
that were first written down by Chai~\cite{Ch}. The use of abelic
rather than abelian varieties also seems to be essential in order to
obtain a good family: this is rather more apparent over a
non-algebraically closed field, when the difference between an abelian
variety (which has a point) and an abelic variety is considerable.

Nakamura, in \cite{Nak2}, takes a different approach. He considers
degenerating families of
abelian varieties with certain types of level structure. In his case
the boundary points correspond to {\em projectively stable
quasi-abelian schemes\/} in the sense of GIT. His construction works over
$\Spec\ZZ[\zeta_N,1/N]$ for a suitable~$N$. At the time of writing it
is not clear whether Nakamura's compactification also leads to the
second Voronoi compactification.

%
%

\vspace{1cm}

\noindent
Authors' addresses:

\bigskip
\noindent
\begin{tabular}{@{}ll@{}}
Klaus Hulek & G. K. Sankaran \\
Institut f\"ur Mathematik & Department of Mathematical Sciences\\
Universit\"at Hannover & University of Bath\\
D 30060 Hannover & Bath BA2 7AY\\
Germany & England\\
{\tt hulek@math.uni-hannover.de} & {\tt gks@maths.bath.ac.uk}
\end{tabular}
\end{document}